\pdfoutput=1
\RequirePackage{ifpdf}
\ifpdf 
\documentclass[pdftex]{sigma}
\else
\documentclass{sigma}
\fi

\usepackage{amsxtra}

\newcommand{\levi}{\nabla^{\mathrm{LC}}}
\DeclareMathOperator{\cocurvature}{cocurv}

\DeclareMathOperator{\affine}{Af\/f}

\DeclareMathOperator{\Adjoint}{Ad}

\DeclareMathOperator{\diff}{Dif\/f}

\DeclareMathOperator{\trace}{trace}

\DeclareMathOperator{\dimension}{dim}

\DeclareMathOperator{\automorphism}{Aut}
\DeclareMathOperator{\selfmorphism}{End}

\DeclareMathOperator{\curvature}{curv}
\DeclareMathOperator{\torsion}{tor}

\newcommand\monomorphism{\hookrightarrow}

\newcommand\identity{\mathrm{id}}
\newcommand\suchthat{\,\left\vert\right.\, }

\renewcommand{\le}{\leqslant}

\numberwithin{equation}{section}

\newtheorem{TheoremBlaom}{Theorem}[section]
\newtheorem{CorollaryBlaom}[TheoremBlaom]{Corollary}
\newtheorem{LemmaBlaom}[TheoremBlaom]{Lemma}
\newtheorem{PropositionBlaom}[TheoremBlaom]{Proposition}
 { \theoremstyle{definition}
\newtheorem{CounterexampleBlaom}[TheoremBlaom]{Counterexample}
\newtheorem{RemarkBlaom}[TheoremBlaom]{Remark} }

\begin{document}

\allowdisplaybreaks

\renewcommand{\PaperNumber}{074}

\FirstPageHeading

\ShortArticleName{The Inf\/initesimalization and Reconstruction of Locally Homogeneous
  Manifolds}

\ArticleName{The Inf\/initesimalization and Reconstruction\\ of Locally Homogeneous Manifolds}

\Author{Anthony D.~BLAOM}

\AuthorNameForHeading{A.D.~Blaom}

\Address{22 Ridge Road, Waiheke Island, New Zealand}
\Email{\href{mailto:anthony.blaom@gmail.com}{anthony.blaom@gmail.com}}

\ArticleDates{Received May 08, 2013, in f\/inal form November 19, 2013; Published online November 26, 2013}

\Abstract{A linear connection on a Lie algebroid is called a
    Cartan connection if it is suitab\-ly compatible with the Lie
  algebroid structure.  Here we show that a smooth connected manifold~$M$ is locally homogeneous~-- i.e., admits an atlas of charts
  modeled on some homogeneous space~$G/H$~-- if and only if there
  exists a transitive Lie algebroid over~$M$ admitting a f\/lat Cartan
  connection that is `geometrically closed'. It is shown how the
  torsion and monodromy of the connection determine the particular
  form of~$G/H$.  Under an additional completeness hypothesis, local
  homogeneity becomes {\em global} homogeneity, up to cover.}


\Classification{53C30;  53C15; 53C07}

\section{Introduction}
Let $M$ be a smooth connected manifold. Then $M$ is {\it locally
  homogeneous} if, for some homogeneous space $G/H$, the smooth
structure of $M$ can be stif\/fened to a ${\mathcal G}$-structure, where
${\mathcal G}$ is the pseudogroup of all those local transformations
of $G/H$ that are restrictions of a left translation by an element of~$G$.

\subsection{Main results}\label{goose}%
The chief purpose of this article is to re-examine local homogeneity
from the Lie groupoid point of view. This leads, in particular, to the
conclusion that a locally homogeneous manifold can be
inf\/initesimalized to obtain a transitive Lie algebroid over $M$,
equipped with a f\/lat linear connection $\nabla$ that is compatible
with the Lie algebroid structure, i.e.\ is a f\/lat {\it Cartan}
connection (see Section~\ref{skoj} below). As we shall elucidate, $\nabla$
being f\/lat amounts to the existence of a~transitive action by some Lie
{\em algebra} on~$M$ `twisted' by a monodromy representation.

Not all transitive Lie algebroids equipped with a f\/lat Cartan
connection are inf\/initesimali\-za\-tions of a locally homogeneous
manifold. If $({\mathfrak g}, \nabla)$ is an inf\/initesimalization,
then $\nabla $ must be {\it geometrically closed}, in a sense made
precise below. Fortunately, geometric closure is suf\/f\/icient for
reversing the inf\/initesimalization procedure:

\begin{TheoremBlaom}\label{TheoremA}
  A smooth manifold $M$ is locally homogeneous if and only if there
  exists a transitive Lie algebroid $\mathfrak g$ over $M$ admitting a
  flat, geometrically closed, Cartan connection $\nabla$.
\end{TheoremBlaom}

As we show in Section~\ref{main}, the particular model $G/H$ that applies is
encoded in the torsion and monodromy of $\nabla$. By contrast, in the
predominant approach to local homogeneity, one f\/ixes a particular
homogeneous space $G/H$ {\em a priori}, and asks if $M$ admits $G/H$
as a local model. In practice, this requires one to anticipate an
appropriate model, or check several candidates systematically. See,
e.g., the survey~\cite{Goldman_10} for this point of view.

\subsection{Geometric closure}\label{geomclose}

Let ${\mathfrak g} $ be a transitive Lie algebroid equipped with a
f\/lat Cartan connection $\nabla$. Fix an arbitrary point $m_0 \in M$
and an arbitrary simply-connected open neighbourhood $U$ of $m_0$. Let
${\mathfrak g}_0$ be the f\/inite-dimensional vector space of $\nabla
$-parallel sections over $U$. Because $\nabla$ is f\/lat, ${\mathfrak
  g}_0$ is the same for all choices of $U$, up to obvious
identif\/ications. Because $\nabla$ is Cartan, ${\mathfrak g}_0 \subset
\Gamma({\mathfrak g}_U)$ is a~subalgebra (with bracket encoded in the
torsion of $\nabla $; see Section~\ref{localform}).

Denote the simply-connected Lie group having ${\mathfrak g}_0 $ as its
Lie algebra by $G_0$, and let ${\mathfrak h}_0$ denote the kernel of
the map
\begin{gather*}
  \xi \mapsto \#\xi (m_0) \colon \  {\mathfrak g}_0 \rightarrow
  T_{m_0}M.
\end{gather*}
Here $\# \colon {\mathfrak g} \rightarrow TM$ denotes the anchor of
$\mathfrak g$. It is easy to see that ${\mathfrak h}_0 \subset
{\mathfrak g}_0$ is a subalgebra, and we say $\nabla $ is {\it
  geometrically closed} if the connected subgroup $H_0 \subset G_0$
with Lie algebra ${\mathfrak h}_0$ is closed in $G_0$. In making this
def\/inition, the choice of f\/ixed point $m_0 \in M$ is immaterial, as we
will establish later in Section~\ref{hudo}.

The basic prototype of a Lie algebroid supporting a f\/lat Cartan
connection is the action algebroid ${\mathfrak g}_0 \times M$
associated with some inf\/initesimal action of a Lie {\em algebra}
${\mathfrak g}_0$ on $M$ (see Section~\ref{morph2}). Indeed, locally this is
the only example (see Section~\ref{localform}). In this case, ${\mathfrak
  h}_0 \subset {\mathfrak g}_0$ is the isotropy subalgebra at $m_0$.

For an example a Lie algebra action that is {\em not} geometrically
closed, see \cite[Example~8]{Kamber_Michor_04} and~\cite{Kowalski_97}.

\subsection{Discussion}
The novelty of Theorem~\ref{TheoremA} lies mainly in the point of view, as
explained in Section~\ref{goose} above. The implications of this change in
viewpoint, as it applies to other parts of Cartan's generalization of
the Klein Erlangen program, are explored in~\cite{Blaom_12}. In its
contemporary conception, this program is described in, e.g.,~\cite{Sharpe_97}.

A Lie algebroid over $M$ equipped with a f\/lat Cartan connection
amounts to a Lie algebra action on the universal cover~$\tilde M$ that
suitably respects covering transformations. We reconstruct a locally
homogeneous structure on $M$ by applying Cartan's development
technique to such actions. Excellent expositions of this technique may
be found elsewhere; see, e.g.,~\cite{Sharpe_97}. That said, Dazord's
integrability result~\cite{Dazord_97}, and Lie groupoid formalism,
provide for an economic treatment of development, which is of\/fered in
Section~\ref{smoothact}.  For a detailed treatment of inf\/initesimal actions
of Lie algebras, we refer the reader to~\cite{Alekseevsky_Michor_95b,Kamber_Michor_04}.

It should be noted that the Lie group $G$ occurring in Theorem~\ref{TheoremA} (and
in Theorem~\ref{TheoremB} mentioned below) need not be connected. Actually, one may
insist that~$G$ be connected (indeed simply-connected) but only at the
cost of allowing transition functions more general than left
translations: rather they may be arbitrary {\em affine}
transformations (for a def\/inition, see Section~\ref{affine}). Our proof of
Theorem~\ref{TheoremA} begins with the proof of a variant along these lines.

Theorem~\ref{TheoremA} suggests a two-step strategy for establishing the local
homogeneity of a smooth manifold: (i)~construct a transitive Lie
algebroid ${\mathfrak g} $ over~$M$, equipped with some Cartan
connection $\nabla$ (not necessarily f\/lat); and (ii)~attempt to deform
$\nabla$, within the class of Cartan connections, to one that is
simultaneously f\/lat and geometrically closed.

For example, associated with any Riemann surface $M$ is a canonical
Lie algebroid ${\mathfrak g} \subset J^1 ({TM})$ (the `isotropy' of
the complex structure; see \cite{Blaom_12}). Associated with any
compatible metric $\sigma$ is a~subalgebroid ${\mathfrak g}_\sigma
\subset {\mathfrak g}$. Applying Cartan's method of equivalence, as we
have described in~\cite{Blaom_12}, one constructs a Cartan connection
$\nabla$ on ${\mathfrak g}_\sigma$, which extends rather naturally to
one on ${\mathfrak g}$.  Applying step~(ii) above to $\nabla$, one can
establish the existence of an atlas of af\/f\/ine coordinates when $M$ has
vanishing Euler characteristic $\chi$; applying (ii) to a
`prolongation' $\nabla^{(1)}$ of~$\nabla $~-- a connection on the
prolongation ${\mathfrak g}^{(1)} \subset J^2({TM})$ of ${\mathfrak
  g}$~-- one obtains an atlas of complex projective coordinates,
provided $\chi > 1$. That is, we recover a weak version of the
uniformization theorem, as in~\cite{Gunning_67}. Details will appear
elsewhere.

In any case, Theorem~\ref{TheoremA} furnishes a strategy for constructing
topological invariants of smooth manifolds, arising as obstructions to
the existence of locally homogeneous structures. Theorem~\ref{TheoremA} also
engenders the following question: What are the global analogues of
{\em intransitive} Lie algebroids over~$M$ equipped with f\/lat Cartan
connections? The answer: A~class of dif\/ferentiable pseudogroups of
transformations on~$M$, generalizing the canonical pseudogroups of
transformations associated with locally homogeneous structures. As it
turns out, such pseudogroups are not necessarily {\em Lie}
pseudogroups in the classical sense. They nevertheless have a very
satisfactory `Lie theory' which is sketched in \cite[Appendix~A]{Blaom_12} and described further in~\cite{Blaom_C}.

\subsection{Completeness and homogeneity up to cover}\label{coco}

Suppose there exists a Lie group $G$ acting transitively on the
universal cover $\tilde M$ of $M$, with $G$ understood to contain the
group $\Gamma\cong \pi_1(M)$ of covering transformations as a
subgroup, and with the action of $G$ on $\tilde M$ extending the
tautological action of~$\Gamma$. In this case $M$ is a double quotient
of groups, $M\cong (G/H)/\Gamma$, and $M$ may be said to be {\it
  homogeneous up to cover}. Under a suitable completeness hypothesis,
locally homogeneous manifolds are already homogeneous up to
cover. See, for example, Thurston's lucid account~\cite{Thurston_97}.

To formulate a notion of completeness for Cartan connections leading
to a strengthening of Theorem~\ref{TheoremA}, let ${\mathfrak g}$ be any Lie
algebroid over $M$ and let $t \mapsto X_t \in {\mathfrak g}$ be a
smooth path, def\/ined on some interval of the real line. Let $m_t\in M$
be the footprint of $X_t$ and $\dot m_t \in TM$ its velocity. Then $t
\mapsto X_t$ is called a $\mathfrak g$-{\it path} if $\#(X_t)=\dot
m_t$, for all $t$. If $\nabla $ is a linear connection on $\mathfrak
g$ (not necessarily f\/lat or Cartan), then a $\mathfrak g$-path $X_t\in
{\mathfrak g}$ is called a {\it geodesic} of $\nabla$ if $\nabla_{\dot
  m_t}X_t=0$ for all $t$. A geodesic of the Levi-Cevita connection
associated with a Riemannian metric is then a geodesic in the standard
sense if we take ${\mathfrak g} = TM$.

The usual argument for the existence of geodesics in Riemannian
geometry carries over to the general case: Through every point of
${\mathfrak g} $ there passes a unique geodesic of $\nabla$. We call
$\nabla$ {\it complete} if every geodesic $X_t$ of $\nabla $ can be
def\/ined for all time $t \in {\mathbb R}$.

The reader is to be warned that compactness of $M$ is {\em not}
suf\/f\/icient for completeness of $\nabla$; a simple counterexample is
given in Section~\ref{recharacterize}. If, in addition, the image of the
monodromy representation has compact closure, then $\nabla $ is indeed
complete. Alternatively, if $M$ admits a complete Riemannian metric
invariant with respect to a natural representation of the Lie
algebroid ${\mathfrak g}$ on $S^2(T^*\!M)$, then $\nabla$ is again
complete. Precise statements and proofs are given in
Section~\ref{completenessSection}, along with a proof of the following
variation of Theorem~\ref{TheoremA}:

\begin{TheoremBlaom}\label{TheoremB}
  A smooth connected manifold $M$ is homogeneous up to cover if and
  only if there exists a transitive Lie algebroid $\mathfrak g$ over~$M$ admitting a flat,  complete, Cartan
  connection~$\nabla$.
\end{TheoremBlaom}

\subsection{Paper outline}

The present article is organized as
follows: In Section~\ref{section1} we review the notion of Cartan
connections on Lie algebroids and explain the sense in which f\/lat
Cartan connections amount to Lie algebra actions twisted by a
monodromy representation (Theorem~\ref{Theorem2.6}). In Section~\ref{sectionI},
we construct the inf\/initesimalization $({\mathfrak g},\nabla)$ of a
locally homogeneous structure on $M$ and observe that it is
geometrically closed. This establishes the necessity of the conditions
in Theorem~\ref{TheoremA}.

In Section~\ref{smoothact} we def\/ine the development of the inf\/initesimal
action of a Lie algebra and describe its behavior under `equivariant
coordinate changes'. This is applied in Section~\ref{main} to reconstruct
a locally homogeneous structure from any transitive Lie algebroid
equipped with a f\/lat, geometrically closed, Cartan connection, which
completes the proof of Theorem~\ref{TheoremA}. In Section~\ref{main} we also explain how the
torsion and monodromy of $\nabla$ determine the particular homogeneous
model~$G/H$ that applies.

In Section~\ref{completenessSection}, we prove Theorem~\ref{TheoremB} and of\/fer some
suf\/f\/icient conditions for completeness. The last section,
Section~\ref{illustration}, illustrates our results by characterising
complete local Lie groups, and by recovering a variant of the
well-known classif\/ication theorem for complete Riemann manifolds.

\subsection{Notation}

Throughout this paper ${\mathfrak g} $ denotes a
Lie algebroid, ${\mathfrak g}_0 $ and ${\mathfrak h}_0$ Lie algebras,
$G_0$, $H_0$, $G$ and $H$ Lie groups, and ${\mathcal G}$ a Lie
groupoid.

\section{Cartan connections on Lie algebroids}\label{section1}

We assume the reader is familiar with the notion of $\mathfrak
g$-connections and $\mathfrak g$-representations ($\mathfrak g$ a Lie
algebroid).
See, for example, \cite{Blaom_05} or~\cite{Blaom_12}.  With the
exception of Section~\ref{holiday}, the present section is mostly a summary
of \cite{Blaom_05}.

\subsection{Cartan connections def\/ined}\label{skoj}

Let ${\mathfrak g} $ be a vector bundle over $M$. Then there is a
one-to-one correspondence between linear connections $\nabla$ on
${\mathfrak g}$ and splittings $s_\nabla \colon {\mathfrak g}
\rightarrow J^1{\mathfrak g} $ of the associated exact sequence
\begin{gather*}
                 0\rightarrow T^*\!M\otimes{\mathfrak g}\monomorphism J^1{\mathfrak
                 g}\rightarrow{\mathfrak g}\rightarrow 0;
\end{gather*}
this correspondence is given by
\begin{gather*}
  s_\nabla X=J^1X + \nabla X.
\end{gather*}
Here $J^1 {\mathfrak g} $ is the vector bundle of one-jets of sections
of $\mathfrak g$, and the inclusion $T^*\!M\otimes{\mathfrak
  g}\monomorphism J^1{\mathfrak g}$ is the morphism whose
corresponding map on section spaces sends $df \otimes X$ to $ f J^1 X
- J^1 (fX)$.

Now suppose that $\mathfrak g$ is not just a vector bundle but a Lie
algebroid. Then $\nabla$ is a {\it Cartan connection} if the vector
bundle morphism $s_\nabla \colon {\mathfrak g} \rightarrow
J^1{\mathfrak g}$ is a morphism {\em of Lie algebroids}.

Recall here that $J^1 {\mathfrak g}$ has a natural Lie algebroid
structure determined by the requirement
\begin{gather*}
  J^1[X,Y]=\big[J^1 X,J^1 Y\big],\qquad \# J^1 X = \# X, \qquad X,Y\in\Gamma({\mathfrak g}).
\end{gather*}
Here and throughout, the anchor of a Lie algebroid is denoted
$\#$. For details and an explicit formula for the bracket on $J^1
{\mathfrak g} $, see~\cite{Blaom_05}, where it is also shown that
$\nabla$ is Cartan if and only if its {\it cocurvature} vanishes. The
latter is a tensor $\cocurvature \nabla \in \Gamma(\wedge^2({\mathfrak
  g}^*)\otimes T^*\!M \otimes {\mathfrak g} )$ def\/ined by
\begin{gather*}
    \cocurvature\nabla ( X, Y)V=\nabla_V[X, Y]-[\nabla_VX,
    Y]-[X,\nabla_VY] +\nabla_{\bar\nabla_XV}Y-\nabla_{\bar\nabla_YV}X.
\end{gather*}

In the above formula $\bar\nabla$ denotes the so-called {\it
  associated ${\mathfrak g}$-connection} on $TM$, def\/ined by
$\bar\nabla_XV=\#\nabla_VX +[\#X,V]$. There is also an {\it associated
  ${\mathfrak g}$-connection} on $\mathfrak g$ itself, also denoted
$\bar\nabla$, and def\/ined by $\bar\nabla_XY=\nabla_YX+[X,Y]$. For this
connection one can def\/ine {\it torsion} in the usual way, by
\begin{gather*}
  \torsion \bar\nabla (X,Y) =
  \bar\nabla_XY-\bar\nabla_YX-[X,Y]=\nabla_{\#Y}X-\nabla_{\#X}Y+[X,Y].
\end{gather*}
When $\nabla $ is Cartan both associated connections are {\em flat},
i.e., def\/ine representations of the Lie algebroid~$\mathfrak g$ on
$TM$ and~${\mathfrak g} $.

\subsection{Cartan connection-preserving morphisms}\label{preserve}

From a well-known characterization of Lie algebroid morphisms given
in, e.g., \cite[Proposition~4.3.12]{Mackenzie_05}, one readily
establishes the following:

\begin{PropositionBlaom}\label{Proposition2.1}
  Let ${\mathfrak g}_1$ be a Lie algebroid over $M_1$ with Cartan
  connection $\nabla^1$ and ${\mathfrak g}_2$ a Lie algebroid over
  $M_2$ with Cartan connection $\nabla^2$. Then a
  connection-preserving vector bundle map $\Phi \colon {\mathfrak g}_1
  \rightarrow {\mathfrak g}_2 $, covering some smooth map $\phi \colon
  M_1 \rightarrow M_2$, is a Lie algebroid morphism if and only~if:
  \begin{enumerate}\itemsep=0pt
  \item[$(1)$] $\# \circ \Phi =T\phi \circ \# $, i.e., $\Phi$ respects
    anchors; and
  \item[$(2)$] $\Phi \torsion \bar \nabla^1 (X,Y)=\torsion \bar\nabla^2(\Phi
    X, \Phi Y) $, for all $X,Y\in{\mathfrak g}_1$.
  \end{enumerate}
  Here $\bar \nabla ^i$ is the associated ${\mathfrak g}_i$-connection
  on ${\mathfrak g}_i$ $($see Section~{\rm \ref{skoj})}.
\end{PropositionBlaom}

\subsection{Equivariance with twist}\label{morph2}

Let ${\mathfrak g}_0$ be a f\/inite-dimensional Lie algebra. If
${\mathfrak g}_0$ acts smoothly on $M$ from the left\footnote{Our
  convention for def\/ining the bracket on the Lie algebra/algebroid of
  a Lie group/groupoid is to use {\em right}-invariant vector f\/ields.}
then we denote the corresponding Lie algebra homomorphism ${\mathfrak
  g}_0 \rightarrow \Gamma(TM)$ by $\xi \mapsto \xi^\dagger$. The
canonical f\/lat connection $\nabla$ on the action algebroid ${\mathfrak
  g}={\mathfrak g}_0 \times M$ is an example of a Cartan
connection. As we recall in Section~\ref{localform} below, this is, locally, the
only example of a f\/lat Cartan connection.

Recall that the anchor of an action algebroid ${\mathfrak g}_0 \times
M$ is def\/ined by $\#(\xi,m)=\xi^\dagger(m)$, and that the bracket on
sections of ${\mathfrak g}_0 \times M$ (${\mathfrak g}_0$-valued
functions on $M$) is given by
\begin{gather}
  [ X, Y]:=\nabla_{\#X}Y-\nabla_{\#Y}X+\tau (X,Y),\label{qwq}
\end{gather}
where $\tau$ is the section of $\wedge^2({\mathfrak g}^*)\otimes
{\mathfrak g} $ def\/ined by
\begin{gather*}
        \tau (X, Y)(m):=[X(m), Y(m)]_{{\mathfrak g}_0},\qquad
        X,Y\in\Gamma(\mathfrak g_0\times M).
\end{gather*}
Note that the associated ${\mathfrak g}$-connection $\bar\nabla$ on
${\mathfrak g}$ has torsion $\tau$.

Now suppose that ${\mathfrak g}_0 $ acts smoothly on two manifolds
$M_1$ and $M_2$, and let $\selfmorphism({\mathfrak g}_0)$ denote the
vector space of Lie algebra endomorphisms of ${\mathfrak g}_0$. Then
we will say that a smooth map $\phi \colon M_1 \rightarrow M_2$ is
{\it ${\mathfrak g}_0$-equivariant} with {\it twist
  $\mu\in\selfmorphism({\mathfrak g}_0)$} if $\xi^\dagger$ and $(\mu
\xi)^\dagger$ are $\phi $-related, for all $\xi\in {\mathfrak g}_0
$. The twist need not be unique.  Applying Proposition~\ref{Proposition2.1}, we
obtain:

\begin{PropositionBlaom}\label{Proposition2.2}
  Every connection-preserving vector bundle morphism
  \[
  {\mathfrak g}_0
  \times M_1 \xrightarrow{\Phi} {\mathfrak g}_0 \times M_2
  \]
 is a Lie
  algebroid morphism if and only if it is of the form $\Phi=\mu \times
  \phi$, for some smooth ${\mathfrak g}_0$-equivariant map $\phi
  \colon M_1 \rightarrow M_2$ with twist $\mu \in
  \selfmorphism({\mathfrak g}_0)$. Here $(\mu \times \phi)(\xi,m):=(\mu \xi,\phi(m))$.
\end{PropositionBlaom}

The signif\/icance of ${\mathfrak g}_0$-equivariance in the present
context is established in Theorem~\ref{Theorem2.6} below.

\subsection{Af\/f\/ine transformations}\label{affine}
Let ${\mathfrak g}_0$ act on $M$ as above, and let $G_0$ be the
simply-connected Lie group integrating ${\mathfrak g}_0$. As we show
in Section~\ref{yut}, if the action of ${\mathfrak g}_0$ is transitive and
geometrically closed, then a ${\mathfrak g}_0$-equivariant map $\phi
\colon M \rightarrow M$ with twist $\mu$ can be `developed' to a map
$\phi_{G_0/H_0} \colon G_0/H_0 \rightarrow G_0/H_0$ on a space of left
cosets $G_0/H_0$. Moreover, $\phi_{G_0/H_0}$ is $G_0$-equivariant with
twist in an appropriate sense. We pause now to def\/ine and characterise
such maps.

For every Lie group $G_0$ one has the group $\affine(G_0)\subset
\diff(M)$ of {\it affine transformations}, generated by left
translations, right translations, and group automorphisms. (Af\/f\/ine
transformations of the Abelian group ${\mathbb R}^n$ are then af\/f\/ine
transformations of~${\mathbb R}^n $ in the standard sense of the
term.) If $H_0 \subset G_0$ is a subgroup, then some elements of~$\affine(G_0)$ descend to bijections of~$G_0/H_0$ that we also refer
to as af\/f\/ine transformations, forming a group denoted by
$\affine(G_0/H_0)$. All left translations in~$G_0$ descend to elements
of $\affine(G_0/H_0)$; a right translation by $k \in G_0$ descends if
and only if $k$ is in the normaliser of~$H_0$; more generally, an
element of $\affine(G_0/H_0)$ is of the form $g \mapsto k\Psi(g)$ for
some $k \in G_0$ and $\Psi \in \automorphism(G_0)$, and descends to an
element of $\automorphism(G_0/H_0)$ if and only if $\Psi(H_0) \subset
H_0$. For example, $\affine({\rm SO}(3)/{\rm SO}(2))\cong
{\rm O}(3)$.

On the other hand, we say that a smooth map $\phi \colon G_0/H_0
\rightarrow G_0/H_0$ is {\em $G_0$-equivariant with twist} $\Psi$ if
$\Psi \colon G_0 \rightarrow G_0$ is a group homomorphism and
  \begin{gather*}
    \phi(g \cdot  x)=\Psi(g) \cdot \phi(x), \qquad g \in G_0,\qquad x \in G_0/H_0.
  \end{gather*}
  The magnanimous reader will readily verify the following
  characterization:
  \begin{PropositionBlaom}\label{Proposition2.3}
  Let $G_0$ and $H_0 \subset G_0$ be connected Lie groups. Then a
  bijection $\phi \colon G_0/H_0 \rightarrow G_0/H_0$ is
  $G_0$-equivariant, with some twist $\Psi \in \automorphism(G_0)$, if
  and only if it is affine.
\end{PropositionBlaom}

\subsection{The local form of a Lie algebroid with f\/lat Cartan
  connection}\label{localform}

For the moment, suppose that $\mathfrak g$ is an arbitrary vector
bundle over~$M$, equipped with a linear connection $\nabla$. Denote
the f\/inite-dimensional subspace of $\nabla $-parallel sections of
$\mathfrak g$ by ${\mathfrak g}_0$. In the special case that $\nabla$
is f\/lat, and $M$ is simply-connected, we obtain a
connection-preserving isomorphism
\begin{gather}
  \begin{split}
 &   {\mathfrak g}_0 \times M\xrightarrow{\sim}{\mathfrak g}, \\
 &   (\xi,m)\mapsto \xi(m).
  \end{split}
  \label{doodah}
\end{gather}
\begin{PropositionBlaom} \label{Proposition2.5}
  Let $\mathfrak g$ be a Lie algebroid over a smooth connected
  manifold $M$ and $\nabla$ a Cartan connection, not necessarily
  flat. Then:
  \begin{enumerate}\itemsep=0pt
  \item[$(1)$] The subspace ${\mathfrak g}_0 \subset \Gamma({\mathfrak g})$
    of $\nabla$-parallel sections is a Lie subalgebra.
  \item[$(2)$] 
  The bracket on ${\mathfrak g}_0$ coincides with
    the torsion of the associated $\mathfrak g$-connection on
    ${\mathfrak g} $, in the sense that
    \begin{gather*}
      [\xi,\eta](m)=\torsion \bar\nabla (\xi(m),\eta(m)),
    \end{gather*}
    for any $m\in M$; $\xi,\eta \in {\mathfrak g}_0$. Here
    $\bar\nabla$ denotes the associated ${\mathfrak g} $-connection on
    $\mathfrak g$.
  \item[$(3)$] 
  The mapping $(\xi,m) \mapsto \# \xi(m) \colon
    {\mathfrak g}_0 \times M \rightarrow TM$ defines a smooth action
    of the Lie algebra~${\mathfrak g}_0$ on~$M$, making ${\mathfrak
      g}_0 \times M$ into an action algebroid.

  \item[$(4)$]
  If $\nabla$ is flat, and $M$ is simply-connected,
    then the canonical isomorphism ${\mathfrak g} \cong {\mathfrak
      g}_0 \times M$ in~\eqref{doodah} is an isomorphism of Lie
    algebroids.
\end{enumerate}
\end{PropositionBlaom}

\subsection{Monodromy and  the global form}\label{holiday}

Again suppose that $\mathfrak g$ is an arbitrary vector bundle over
$M$. Let $\tilde {\mathfrak g} $ denote its pullback to a vector
bundle over the universal cover~$\tilde M$ of~$M$. Let $\Gamma$ denote
the group of covering transformations of~$M$. Then it is an elementary
observation that there exists a unique lift of the tautological action
of $\Gamma$ on $\tilde M$ to an action on $\tilde {\mathfrak g} $
satisfying the following requirement: For all $X\in \tilde {\mathfrak
  g} $ and $\phi \in \Gamma$, $X$~and~$\phi \cdot X $ have the same
image under the canonical projection $\tilde {\mathfrak g} \rightarrow
{\mathfrak g}$. Evidently $\Gamma$ acts on~$\tilde {\mathfrak g} $ by
vector bundle automorphisms, and we recover the original vector bundle
over~$M$ as a quotient: ${\mathfrak g}=\tilde {\mathfrak g}/\Gamma $.

Now suppose $\mathfrak g$ is equipped with a f\/lat linear connection~$\nabla$ and let~$\tilde \nabla $ denote the pullback of~$\nabla $ to
a connection on $\tilde {\mathfrak g} $, also f\/lat. Let ${\mathfrak
  g}_0 \subset \Gamma(\tilde {\mathfrak g})$ denote the vector space
of $\tilde \nabla $-parallel sections. Then, as $\tilde M$ is
simply-connected, we have a canonical isomorphism $\tilde {\mathfrak
  g} \cong {\mathfrak g}_0 \times \tilde M$ (as discussed in
Section~\ref{localform} above) and consequently ${\mathfrak g} \cong
({\mathfrak g}_0 \times \tilde M)/\Gamma$.  Since the action of
$\Gamma $ on $\tilde {\mathfrak g} $ described above automatically
preserves $\tilde \nabla $, it must be of the form
\begin{gather}
  \phi \cdot (\xi, \tilde m)=(\mu_\phi \xi,\phi(\tilde m))=(\mu_\phi \times \phi)(\xi,\tilde m),\label{form}
\end{gather}
for some uniquely determined group homomorphism $\mu \mapsto \mu_\phi
\colon \Gamma \rightarrow \automorphism({\mathfrak g}_0)$. This is the
{\it mo\-nodromy representation} associated with the f\/lat connection $\nabla$.

\begin{RemarkBlaom}\label{Remark2.6}
  As the reader will recall, monodromy has the following alternative
  interpretation. Let $\tilde m_0 \in \tilde M$ be a point covering
  any f\/ixed point $m_0 \in M$. Then there is an isomorphism
  ${\mathfrak g}_0 \cong {\mathfrak g}|_{m_0}$ in which each $\xi \in
  {\mathfrak g}_0 $ corresponds to the image of $\xi(\tilde m_0)$
  under the canonical projection $\tilde {\mathfrak g} \rightarrow
  {\mathfrak g}$. Also, $\Gamma$ may be identif\/ied with the
  fundamental group $\pi_1(M,m_0)$. Under these identif\/ications~$\mu_\phi (\xi )$ is the $\nabla$-parallel translate of $\xi \in
  {\mathfrak g}|_{m_0}$ along any closed path $\gamma $ in $M$
  representing $\phi \in \pi_1(M, m_0)$. Or, this parallel translation
  may be viewed in the following way: each $\xi \in {\mathfrak
    g}|_{m_0}$ extends to a {\em locally defined} $\nabla$-parallel
  section $X$ of $\mathfrak g$, which can be `analytically
  continued' around $\gamma $ and re-evaluated at~$m_0$ to obtain
  $\mu_\phi (\xi)$.
\end{RemarkBlaom}

\begin{PropositionBlaom}\label{Proposition2.6}
  In the scenario above, suppose that $\mathfrak g$ is a Lie algebroid
  and $\nabla$ a flat Cartan connection on $\mathfrak g$, so that
  ${\mathfrak g}_0$ is a Lie algebra acting on $\tilde M$ $($see
  Proposition~{\rm \ref{Proposition2.5})}. Then:
\begin{enumerate}\itemsep=0pt
\item[$(1)$]
The canonical isomorphism ${\mathfrak g} \cong ({\mathfrak g}_0
  \times \tilde M)/\Gamma$ is an isomorphism of Lie algebroids.
\item[$(2)$]
The group of covering transformations $\Gamma$ acts in the
  monodromy representation $\mu \mapsto \mu_\phi$ by Lie algebra
  automorphisms of ${\mathfrak g}_0$.
\item[$(3)$]
Each covering transformation $\phi \in \Gamma $ is a
  ${\mathfrak g}_0 $-equivariant diffeomorphism with twist
  $\mu_\phi\in \automorphism({\mathfrak g}_0)$.
\item[$(4)$]
If we restrict the tensor $\torsion \bar\nabla \in
  \Gamma(\wedge^2({\mathfrak g}^*)\otimes {\mathfrak g})$ to define a
  bracket on the fibre ${\mathfrak g}|_{m_0}$, then the isomorphism
  ${\mathfrak g}_0\cong {\mathfrak g}|_{m_0}$ in Remark~{\rm \ref{Remark2.6}} is
  an isomorphism of Lie algebras. Here $\bar\nabla $ denotes the
  $\mathfrak g$-connection on $\mathfrak g$ associated with $\nabla$.
\end{enumerate}
\end{PropositionBlaom}

\begin{proof}
  Evidently, there is a unique Lie algebroid structure on $\tilde
  {\mathfrak g}$ such that the projection $\tilde {\mathfrak g}
  \rightarrow {\mathfrak g} $ is a Lie algebroid morphism. With
  respect to this structure, the action of $\Gamma$ on $\tilde
  {\mathfrak g} $ is by Lie algebroid isomorphisms. Thus the
  isomorphism ${\mathfrak g} \cong \tilde {\mathfrak g} / \Gamma $ may
  be regarded as a Lie algebroid isomorphism. (None of these
  statements depend on the existence of a trivialization $\tilde
  {\mathfrak g} \cong {\mathfrak g}_0 \times \tilde M$.) On the other
  hand, the isomorphism $\tilde {\mathfrak g} \cong {\mathfrak g}_0
  \times \tilde M$ determined by the f\/lat connection $\tilde \nabla $
  is a Lie algebroid isomorphism, by Proposition~\ref{Proposition2.5}(4).
   So~(1) holds.

  Conclusions (2) and (3) are consequences of
  Proposition~\ref{Proposition2.2}, 
  and the fact that $\Gamma $ acts on
  $\tilde {\mathfrak g} \cong {\mathfrak g}_0 \times \tilde M$ by Lie
  algebroid isomorphisms. Conclusion (4) follows from
  Proposition~\ref{Proposition2.5}(2) 
  and Proposition~\ref{Proposition2.1}. 
\end{proof}

The constructions above are reversible. Indeed, let ${\mathfrak g}_0$
be an arbitrary Lie algebra acting smoothly on the universal cover
$\tilde M$ of $M$, and let $\phi \mapsto \mu_\phi \colon \Gamma
\rightarrow \automorphism ({\mathfrak g}_0)$ be a representation of~$\Gamma$ by Lie algebra automorphisms of ${\mathfrak g}_0$ satisfying~(3) 
(of which there may be more than one). Then the action of
$\Gamma$ on the action algebroid $\tilde {\mathfrak g} :={\mathfrak
  g}_0 \times \tilde M$ def\/ined by \eqref{form} is by Lie algebroid
automorphisms (by Proposition~\ref{Proposition2.2}), 
implying that the
quotient ${\mathfrak g}:=\tilde {\mathfrak g} /\Gamma$ is a Lie
algebroid. Moreover, the canonical f\/lat Cartan connection on $\tilde
{\mathfrak g} ={\mathfrak g}_0 \times \tilde M$ drops to a f\/lat Cartan
connection on ${\mathfrak g}$ whose monodromy is precisely~$\mu$. This
establishes the following:

\begin{TheoremBlaom}\label{Theorem2.6}
  Let $\Gamma \cong \pi_1(M)$ denote the group of covering
  transformations of $\tilde M$. Then there is a natural one-to-one
  correspondence between: $(i)$ pairs $({\mathfrak g}, \nabla$), where
  $\mathfrak g$ is a Lie algebroid over~$M$ and~$\nabla$ is a flat
  Cartan connection; and $(ii)$ pairs $({\mathfrak g}_0,\mu)$, where
  ${\mathfrak g}_0$ is a finite-dimensional Lie algebra acting
  smoothly on $\tilde M$, and $\mu$ is a representation of $\Gamma $
  on ${\mathfrak g}_0$ by Lie algebra automorphisms such that each
  covering transformation $\phi \in \Gamma $ is a ${\mathfrak g}_0
  $-equivariant diffeomorphism with twist $\mu_\phi\in
  \automorphism({\mathfrak g}_0)$.
\end{TheoremBlaom}

\section{Inf\/initesimalization}\label{sectionI}
In this section we prove necessity of the conditions in Theorem~\ref{TheoremA},
i.e., for every locally homogeneous manifold $M$, there exists a
transitive Lie algebroid ${\mathfrak g}$ over $M$ supporting a f\/lat,
geometrically closed, Cartan connection $\nabla$. We refer to the
particular pair $({\mathfrak g}, \nabla)$ constructed below as the
locally homogeneous structure's {\it infinitesimalization}.

\subsection{Construction of the
  inf\/initesimalization}\label{infinitesimalization}
Suppose $M$ is a locally homogeneous manifold modeled on $G/H$ and let
$\psi_i \colon U_i \rightarrow G/H$, $i\in I$, be an atlas of
coordinate charts adapted to the model. This means that for each $i,j
\in I$ for which $U_i\cap U_j \ne \varnothing$, we are given an element
$g_{ji}\in G$, such that $\big(\psi_j \circ \psi_i^{-1}\big)(x) = g_{ji}\cdot
x$ for each $x\in \psi_i(U_i)\cap \psi_j(U_j)$. Moreover, the elements
$g_{ji}$ satisfy the cocycle conditions $g_{ii}=\identity$, and
$g_{ij}g_{jk}=g_{ik}$ for all $i,j,k \in I$.

Now let ${\mathfrak g}_0$ denote the Lie algebra of
$G$. Inf\/initesimalizing the action of $G$ on $G/H$, we obtain a Lie
algebra homomorphism $\xi \mapsto \xi^\dagger \colon {\mathfrak g}_0
\rightarrow \Gamma(T(G/H)) $, which makes ${\mathfrak g}_0 \times G/H$
into a transitive action algebroid. In particular, each restriction
${\mathfrak g}_0 \times \psi_i(U_i) \subset {\mathfrak g}_0 \times
G/H$ is a transitive action algebroid over $\psi_i(U_i)$. We claim
that the Lie algebroids ${\mathfrak g}_0 \times \psi_i(U_i) $, $i \in
I$, are local trivializations of a single transitive Lie algebroid
${\mathfrak g} $ over $M$, the canonical f\/lat connections on the
${\mathfrak g}_0 \times \psi_i(U_i)$ representing a f\/lat Cartan
connection $\nabla $ on $\mathfrak g$.

To see this, def\/ine an equivalence relation $\sim$ on the
set
\[
\{(\xi,x,i) \in {\mathfrak g}_0 \times G/H \times I\suchthat x
\in \psi_i(U_i)\}
\] (the disjoint union of the sets ${\mathfrak g}_0
\times \psi_i(U_i)$, $i\in I$) by declaring
$(\xi,x,i)\sim(\xi',x',i')$ whenever there is an $m\in U_i\cap U_{i'}$
such that $\psi_i(m)=x$, $\psi_{i'}(m)=x'$ and
$\xi'=\Adjoint_{g_{i'i}}\xi$. Then the set ${\mathfrak g}$ of
equivalence classes is a smooth vector bundle over $M$, with footprint
projection $[\xi,x,i]\mapsto \psi_i^{-1}(x)$; here $[\xi,x, i]$
denotes the class with representative $(\xi,x,i)$.

The vector bundle $\mathfrak g$ admits local trivializations given by
\begin{gather*}
  {\mathfrak g}|_{U_i} \xrightarrow{\Psi_i}{\mathfrak g}_0 \times \psi_i(U_i), \\
  [\xi,x,i] \mapsto (\xi,x),
\end{gather*}
with transition functions given by
\begin{gather*}
  \big(\Psi_j \circ
  \Psi_i^{-1}\big)(\xi,x)=(\Adjoint_{g_{ji}}\xi,g_{ji}\cdot x).
\end{gather*}
In particular, the transition functions preserve the canonical f\/lat
connections on the action algebroids ${\mathfrak g}_0 \times
\psi_i(U_i)$. It follows that there is a (necessarily f\/lat) connection
$\nabla$ on $\mathfrak g$ that is locally represented by the canonical
f\/lat connection on each ${\mathfrak g}_0 \times \psi_i(U_i)$.

If we write $L_g(x):=g\cdot x$, then the pushforward of $\xi^\dagger$
by the transformation $L_g \colon G/H \rightarrow G/H$ is $(\Adjoint_g
\xi)^\dagger$. Using this fact, one easily sees that there is a
well-def\/ined vector bundle epimorphism $\# \colon {\mathfrak
  g} \rightarrow TM$, def\/ined locally by
\begin{gather*}
  \#\big(\Psi_i^{-1}(\xi,x)\big)=T \psi_i^{-1}\cdot \xi^\dagger(x).
\end{gather*}
It remains to def\/ine a Lie bracket on sections of $\mathfrak g$ for
which $\#$ is a compatible anchor. To this end, notice that each local
trivialization $\Psi_i$ identif\/ies a f\/ibre ${\mathfrak g}|_m$ ($m\in
U_i$) with ${\mathfrak g}_0$. This identif\/ication depends on the local
trivialization chosen, but only up to adjoint transformations of~${\mathfrak g}_0$; it consequently transfers the Lie bracket on
${\mathfrak g}_0$ to one on ${\mathfrak g}|_m$ that is
trivialization-independent. We let $\tau\in\wedge^2({\mathfrak
  g}^*)\otimes {\mathfrak g}$ denote the tensor whose restriction to
each f\/ibre ${\mathfrak g}|_m$ is the Lie bracket just def\/ined. A~bracket on ${\mathfrak g}$ is then given by
(cf.~\eqref{qwq}):
\begin{gather*}
  [ X, Y]:=\nabla_{\#X}Y-\nabla_{\#Y}X+\tau (X,Y).
\end{gather*}
With this bracket $\mathfrak g$ becomes a transitive Lie algebroid and
the local trivializations $\Psi_i \colon {\mathfrak g}|_{U_i}
\rightarrow {\mathfrak g}_0 \times \psi_i(U_i)$ become
connection-preserving Lie algebroid morphisms. In particular, the
connec\-tion~$\nabla$ is Cartan and geometrically closed because its
local representatives are.

\section{The development of Lie algebra actions}\label{smoothact}
To reconstruct a locally homogeneous structure, from a pair
$({\mathfrak g},\nabla)$ satisfying the conditions in Theorem~\ref{TheoremA}, we
will use Cartan's development technique. We pause here to describe
development using the economy af\/forded by Lie groupoid language, and
to show that development is suitably equivariant. Before doing so, we
argue that geometric closure, as def\/ined in Section~\ref{geomclose}, is
independent of the choice of f\/ixed point $m_0 \in M$.

\subsection{Geometric closure recharacterized}\label{hudo}

Evidently a f\/lat Cartan connection $\nabla$ on a transitive Lie
algebroid $\mathfrak g$ over $M$ is geometrically closed `at the point
$m_0 \in M$' if and only if $\tilde \nabla$ is geometrically closed at
some point $\tilde m_0 \in \tilde M$ covering $m_0$. Here $(\tilde
{\mathfrak g},\tilde \nabla )$ denotes the lift of $({\mathfrak
  g},\nabla )$ to the universal cover $\tilde M$, as described in
Section~\ref{holiday}. So, without loss of generality, we now suppose $M$ is
simply-connected.

According to Proposition~\ref{Proposition2.5}(4),
  $\mathfrak g$ is
isomorphic to an action algebroid ${\mathfrak g}_0 \times M$, where
${\mathfrak g}_0$ has the same meaning as in
Section~\ref{geomclose}. According to Dazord~\cite{Dazord_97}, all action
algebroids are integrable. So there is a Lie groupoid ${\mathcal G} $
integrating $ {\mathfrak g}_0 \times M$, which we may take to be
source-simply-connected \cite[Lie~I]{Crainic_Fernandes_11}.

Fix some $m_0 \in M$ and let ${\mathcal G}_{m_0}^{m_0}$ denote the
isotropy group at $m_0$ (the group of arrows of ${\mathcal G}$
simultaneously beginning and ending at~$m_0$).

\begin{LemmaBlaom}\label{Lemma4.1}
The isotropy group ${\mathcal G}_{m_0}^{m_0}$ is connected.
\end{LemmaBlaom}

\begin{proof}
  Since we assume ${\mathfrak g}_0$ acts transitively, ${\mathcal G} $
  is a transitive Lie groupoid (because its orbits are disjoint and
  open and~$M$ is connected). Consequently, if $P \subset {\mathcal
    G}$ denotes the source-f\/ibre over~$m_0$ (the subset of all arrows
  in ${\mathcal G} $ beginning at~$m_0$), then $P$ is a principal
  ${\mathcal G}_{m_0}^{m_0}$-bundle over~$M$, with the Lie group
  ${\mathcal G}_{m_0}^{m_0}$ acting on~$P$ from the right. The bundle
  projection $P \rightarrow M$ is just the restriction of the
  target-projection ${\mathcal G} \rightarrow M$ of the groupoid~$\mathcal G$. For this principal bundle we have a corresponding long
  exact sequence in homotopy,
  \begin{gather*}
    \cdots \rightarrow  \pi_1({\mathcal G}_{m_0}^{m_0})\rightarrow \pi_1(P)
    \rightarrow \pi_1(M) \rightarrow \pi_0 ({\mathcal
      G}_{m_0}^{m_0})\rightarrow \pi_0(P) \rightarrow \cdots  .
  \end{gather*}
  Since $P$ is connected and $M$ is simply-connected, $\pi_0({\mathcal
    G}_{m_0}^{m_0})$ is trivial.
\end{proof}

Let $G_0$ and $H_0$ have the meanings given in Section~\ref{geomclose}. Then
then there exists a Lie groupoid morphism $\Omega \colon {\mathcal G}
\rightarrow G_0$ integrating the canonical projection ${\mathfrak g}_0
\times M \rightarrow {\mathfrak g}_0$ (a morphism of Lie algebroids)
\cite[Lie~II]{Crainic_Fernandes_11}. Evidently the subgroups
$\Omega({\mathcal G}_{m_0}^{m_0})$ and~$H_0$ of~$G_0$ have the same
Lie algebra~${\mathfrak h}_0$. By the lemma they must coincide. This
proves:
\begin{PropositionBlaom}\label{Proposition4.2}
  The Cartan connection $\nabla$ is geometrically closed if and only
  if the $($necessarily connected$)$ Lie group $\Omega({\mathcal
    G}^{m_0}_{m_0})$ is closed in $G_0$.
\end{PropositionBlaom}

The independence of the choice of point $m_0 \in M$ is now clear: If
$m_0' \in M$ is a second point then, by the transitivity of ${\mathcal
  G} $, there exits an arrow $p \in {\mathcal G} $ from $m_0$ to
$m_0'$, in which case ${\mathcal G}^{m_0'}_{m_0'}=p {\mathcal
  G}^{m_0}_{m_0} p^{-1}$ and $\Omega\big({\mathcal G}^{m_0'}_{m_0'}\big)=g
\Omega({\mathcal G}^{m_0}_{m_0})g^{-1}$, where $g=\Omega(p)$.

\subsection{Development def\/ined}\label{devmap}
Let ${\mathfrak g}_0$ be any f\/inite-dimensional Lie algebra acting
smoothly on $M$.  {\em Assume the action is transitive and
  geometrically closed.}  The development is always def\/ined as a map
from the universal cover $\tilde M$ (on which ${\mathfrak g}_0$ acts
also) so, without loss of generality, we suppose once more that $M$ is
simply-connected. We again denote by $G_0$ the simply-connected Lie
group having ${\mathfrak g}_0 $ as Lie algebra, and let ${\mathcal
  G}$, ${\mathcal G}^{m_0}_{m_0}$, $P$, and $\Omega \colon {\mathcal
  G} \rightarrow G_0$ have the meanings given in the preceding
Section~\ref{hudo}. By Proposition~\ref{Proposition4.2} 
and geometric closure,
$H_0:=\Omega({\mathcal G}^{m_0}_{m_0})$ is a closed subgroup of~$G_0$,
so that~$G_0/H_0$ is a smooth Hausdorf\/f manifold.

Because $\Omega \colon {\mathcal G} \rightarrow G_0$ integrates the
projection ${\mathfrak g}_0 \times M \rightarrow {\mathfrak g}_0$ (a
point-wise isomorphism), its restriction to $P$ is a local
dif\/feomorphism $\Omega \colon P\rightarrow G_0$. Using the fact that
$\Omega$ is a groupoid morphism, we see that $\Omega$ sends orbits of
${\mathcal G}^{m_0}_{m_0}$ in $P$ to orbits of $H_0$ in $G_0$ (left
cosets). Since $P/{\mathcal G}^{m_0}_{m_0}\cong M$, it follows that
$\Omega\colon P\rightarrow G_0$ descends to a map $D \colon
M\rightarrow G_0/H_0$:
\begin{gather*}
  \begin{CD}
    P   @>{\Omega}>> G_0\\
    @VV{/{\mathcal G}^{m_0}_{m_0}}V @VV{/H_0}V \\
    M@>{\exists\, D}>> G_0/H_0
  \end{CD}
\end{gather*}
This is the {\it development} determined by the arbitrary choice of
point $m_0$. One has $D(m)=\Omega(p)H_0$, where $p \in {\mathcal G}$
is any arrow from $m_0$ to $m$.
\begin{PropositionBlaom}\label{Proposition4.3}
  The development $D \colon M \rightarrow {G_0/H_0}$ is a local
  diffeomorphism.
\end{PropositionBlaom}

\begin{proof}
  This is a straightforward consequence of the fact that the principal
  bundle projection $P \rightarrow M$ is a surjective submersion (and
  so admits local sections), $\Omega \colon P \rightarrow G_0$ is a
  local dif\/feomorphism, and $\dimension ({\mathcal
    G}^{m_0}_{m_0})=\dimension (H_0)$.
\end{proof}

Restricting the development to suf\/f\/iciently small open sets in $M$, we
obtain an atlas of charts trivially adapted the homogeneous space
${G_0/H_0}$:
\begin{CorollaryBlaom}\label{Corollary4.4}
  If a finite-dimensional Lie algebra ${\mathfrak g}_0$ acts on a
  simply-connected manifold~$M$, and this action is transitive and
  geometrically closed, then~$M$ is locally homogeneous. Indeed in
  that case there exists a closed subgroup~$H_0 \subset G_0$ of the
  simply-connected Lie group integrating~${\mathfrak g}_0$, and an
  atlas of charts adapted to the homogeneous model~$G_0/H_0$ whose
  transition functions are all identity transformations.
\end{CorollaryBlaom}

\subsection{Behavior under equivariant coordinate changes}\label{yut}
In order to describe how development transforms under
equivariant coordinate changes, we wish to associate, with each
${\mathfrak g}_0$-equivariant dif\/feomorphism $\phi \colon M
\rightarrow M$ with twist $\mu \in \automorphism ({\mathfrak g}_0)$, a~corresponding dif\/feomorphism
\[
\phi_{G_0/H_0}\colon \ {G_0/H_0}
\rightarrow {G_0/H_0}.
\]
 To this end, we require the following:

\begin{LemmaBlaom}\label{LemmaA}
  Let $\phi \colon M \rightarrow M$ be any smooth ${\mathfrak
    g}_0$-equivariant map with twist $\mu \in \selfmorphism({\mathfrak
    g}_0)$. Let $\hat\mu \colon G_0 \rightarrow G_0$ denote the unique
  group homomorphism with derivative $\mu$, and let $(\mu \times
  \phi)^\wedge \colon {\mathcal G} \rightarrow {\mathcal G} $ denote
  the unique Lie groupoid morphism with derivative $\mu \times \phi \colon
  {\mathfrak g}_0 \times M \rightarrow {\mathfrak g}_0 \times M$. Then
  the following diagram commutes:
\begin{gather*}
  \begin{CD}
    {\mathcal G}   @>{(\mu \times \phi)^\wedge}>> {\mathcal G}\\
    @VV{\Omega}V @VV{\Omega}V \\
    G_0@>{\hat\mu}>> G_0
  \end{CD}
\end{gather*}
\end{LemmaBlaom}

\begin{proof}
  The Lie groupoid morphisms $\Omega \circ (\mu \times \phi)^\wedge$ and
  $\hat\mu \circ \Omega$ have the same derivatives, namely
  $(\xi,m)\mapsto \mu \xi \colon {\mathfrak g}_0 \times M \rightarrow
  {\mathfrak g}_0 $. By the uniqueness part of the generalization to
  Lie groupoids of Lie's second integrability theorem \cite[Lie~II]{Crainic_Fernandes_11}, these morphisms must coincide.
\end{proof}

\begin{LemmaBlaom}\label{LemmaB}
 Suppose $\phi \colon M \rightarrow M$ is a ${\mathfrak
    g}_0$-equivariant diffeomorphism with twist $\mu \in
  \automorphism({\mathfrak g}_0)$. Let $q \in {\mathcal G} $ be any
  arrow from $m_0$ to $\phi(m_0)$, and let $\hat\mu$ be the unique
  automorphism of $G_0$ with derivative $\mu$. Then
  $\hat\mu(H_0)=\Omega(q)H_0 \Omega(q)^{-1}$.
\end{LemmaBlaom}
\begin{proof}
  Because, in the notation of Lemma~\ref{LemmaA}, the groupoid
  automorphism $(\phi \times \mu)^\wedge \colon {\mathcal G}
  \rightarrow {\mathcal G} $ covers $\phi \colon M \rightarrow M$, and
  because $q \in {\mathcal G} $ is an arrow from $m_0$ to $\phi(m_0)$,
  we have $(\phi \times \mu)^\wedge({\mathcal G}_{m_0}^{m_0})=q
  {\mathcal G}_{m_0}^{m_0}q^{-1}$. Applying that lemma, we
  compute
  \begin{gather*}
    \hat\mu(H_0)=\hat\mu\big(\Omega({\mathcal G}_{m_0}^{m_0})\big)=\Omega\big( (\phi \times \mu)^\wedge({\mathcal G}_{m_0}^{m_0}) \big)=\Omega\big(q {\mathcal G}_{m_0}^{m_0}q^{-1}\big)=\Omega(q)H_0 \Omega(q)^{-1}.\tag*{\qed}
  \end{gather*}
  \renewcommand{\qed}{}
\end{proof}

Lemma~\ref{LemmaB} implies that the map $\phi_{G_0/H_0}\colon {G_0/H_0}
\rightarrow {G_0/H_0}$, given implicitly by
\begin{gather}
  \phi_{G_0/H_0} (g H_0)=\hat\mu(g)\Omega(q)H_0
  \label{abox}
\end{gather}
is well-def\/ined. It is also independent of the choice of arrow $q \in
{\mathcal G} $ from $m_0$ to $\phi(m_0)$, because
$\Omega(q)H_0=D(\phi(m_0))$. Note that despite our choice of notation,
$\phi_{G_0/H_0}$ depends not just on~$\phi $ but also on the twist
$\mu$.

\begin{PropositionBlaom}\label{Proposition4.7} Let $D \colon M \rightarrow {G_0/H_0}$ be the
  development associated with some finite-di\-men\-sio\-nal Lie algebra
  ${\mathfrak g}_0$ acting transitively on $M$. Then, for all
  ${\mathfrak g}_0$-equivariant diffeomorphisms $\phi,\phi^\prime
  \colon M \rightarrow M$ with twists $\mu$, $\mu^\prime$ respectively,
  one has:
  \begin{enumerate}\itemsep=0pt
  \item[$(1)$] 
  $\phi_{G_0/H_0} \circ \phi^\prime_{G_0/H_0}=(\phi \circ
    \phi^\prime)_{G_0/H_0}$.

  \item[$(2)$] 
  $\phi_{G_0/H_0} \colon {G_0/H_0} \rightarrow
    {G_0/H_0} $ is $G_0$-equivariant with twist $\hat\mu$, in the
    sense that
  \begin{gather*}
    \phi_{G_0/H_0}(g\cdot x)=\hat\mu(g)\cdot \phi_{G_0/H_0}(x), \qquad g
    \in G_0,\qquad x \in {G_0/H_0}.
  \end{gather*}

  \item[$(3)$] 
  The following diagram commutes:
    \begin{gather*}
      \begin{CD}
          M   @>{\phi}>> M\\
          @V{D}VV @VV{D}V \\
          {G_0/H_0}@>{\phi_{G_0/H_0}}>> {G_0/H_0}
        \end{CD}
    \end{gather*}
  \end{enumerate}
   Here $\hat \mu$ denotes the unique automorphism of $G_0$
   integrating the Lie algebra morphism $\mu \colon {\mathfrak g}_0
   \rightarrow {\mathfrak g}_0$.
 \end{PropositionBlaom}

 In other words, if $\automorphism({\mathfrak g}_0 \times M)$ denotes
 the group of all ${\mathfrak g}_0$-equivariant dif\/feomorphisms of~$M$
 with twist, then the map $(\phi,\mu) \mapsto \phi_{G_0/H_0}$ def\/ines
 an action of $\automorphism({\mathfrak g}_0 \times M)$ on ${G_0/H_0}$
 by af\/f\/ine transformations (see Proposition~\ref{Proposition2.3}) 
 and with
 respect to this action and the tautological action of
 $\automorphism({\mathfrak g}_0 \times M)$ on $M$, the development $D
 \colon M \rightarrow {G_0/H_0}$ is equivariant, in the standard sense
 of group actions.

\begin{proof}
  Let $q,q^\prime \in {\mathcal G} $ be arrows from $m_0$ to
  $\phi(m_0)$, $\phi^\prime(m_0)$ respectively. Then
  $q^{\prime\prime}:=(\phi \times \mu)^\wedge(q^\prime)\,q$ is an
  arrow from $m_0$ to $\phi(\phi^\prime(m_0))$. It follows that for an
  arbitrary element $gH_0\in {G_0/H_0}$, we have
  \begin{gather*}
    (\phi \circ  \phi^\prime)_{G_0/H_0}(gH_0)=\widehat{\mu \mu^\prime}(g)\Omega(q^{\prime\prime})H_0=\hat\mu\big(\widehat{\mu^\prime}(g)\big)\Omega\big( (\phi \times \mu)^\wedge(q^\prime)q \big)H_0\\
    \hphantom{(\phi \circ  \phi^\prime)_{G_0/H_0}(gH_0)}{}
=\hat\mu\big(\widehat{\mu^\prime}(g)\big)\hat\mu(\Omega(q^\prime))\Omega(q)H_0=
\hat\mu \big( \widehat{\mu^\prime}(g)\Omega(q^\prime) \big)\Omega(q)H_0,
  \end{gather*}
  where at the beginning of the second line we have applied Lemma~\ref{LemmaA}. On the other hand, we have
  \begin{gather*}
    \phi_{G_0/H_0}(\phi^\prime_{G_0/H_0}(gH_0))=
     \phi_{G_0/H_0}(\widehat{\mu^\prime}(g)\Omega(q^\prime)H_0)=\hat\mu \Big(\,\widehat{\mu^\prime}(g)\Omega(q^\prime)\,\Big)\Omega(q)H_0.
   \end{gather*}
   Comparing this equation with the preceding one establishes~(1).

   One deduces (2) 
   immediately from the def\/inition of
   $\phi_{G_0/H_0} $. Regarding~(3), 
   let $m \in M$ be
   arbitrary and let $p \in {\mathcal G} $ be an arrow from $m_0$ to
   $m$. Then $p^\prime := (\phi \times \mu)^\wedge(p) q$ is an arrow
   from $m_0$ to $\phi(m)$. Consequently, we compute
\begin{gather*}
  D(\phi(m))=\Omega(p^\prime)H_0=\Omega \big( (\phi \times
  \mu)^\wedge(p) \big)\Omega(q)H_0\\
\hphantom{D(\phi(m))}{}  =\hat\mu(\Omega(p))\Omega(q)H_0
  =\phi_{G_0/H_0}(\Omega(p)H_0)=\phi_{G_0/H_0}(D(m)).
\end{gather*}
At the beginning of the second line we have again applied Lemma~\ref{LemmaA}.
\end{proof}

\section{Reconstruction}\label{main}

In this section we complete the proof of Theorem~\ref{TheoremA} by reconstructing a
locally homogeneous structure from any transitive Lie algebroid
supporting a f\/lat, geometrically closed Cartan connection. We also
explain how torsion and monodromy of the connection determine a
suitable model~$G/H$.

\subsection{The role of torsion and monodromy}\label{rolltorsion}
According to the following result, if~$M$ satisf\/ies the hypotheses of
Theorem~\ref{TheoremA}, then it locally homogeneous in the apparently weaker sense
of admitting an atlas in which the transition functions are {\em
  affine} transformations of a homogeneous model $G_0/H_0$ (in the
sense of Section~\ref{affine}). The advantage of this formulation over the one
in Theorem~\ref{TheoremA} is that we may take $G_0$ to be simply-connected and
$H_0$ to be connected.

\begin{PropositionBlaom}\label{Proposition5.1}
  Let $\mathfrak g$ be a transitive Lie algebroid, over a smooth
  connected manifold $M$, and~$\nabla $ a~flat Cartan connection on
  $\mathfrak g$. Define a ${\mathfrak g}$-connection $\bar\nabla$ on
  $\mathfrak g$ by
  \begin{gather*}
    \bar\nabla_XY=\nabla_{\#Y}X+[X,Y],
  \end{gather*}
  and let $\torsion \bar\nabla $ denotes its torsion:
  \begin{gather*}
    \torsion \bar\nabla (X,Y)=\bar\nabla_XY-\bar\nabla_YX-[X,Y]=\nabla_{\#Y}X-\nabla_{\#X}Y+[X,Y].
  \end{gather*}
  Then, fixing a point $m_0 \in M$, we have:
  \begin{enumerate}\itemsep=0pt
  \item[$(1)$]  
  The restriction of $\torsion \bar\nabla$ to ${\mathfrak
      g}_0:={\mathfrak g}|_{m_0}$ is a Lie bracket on ${\mathfrak
      g}_0$.
  \end{enumerate}

  Next, let $\Gamma=\pi_1(M,m_0)$ denote the fundamental group of $M$
  and let $\mu \colon \Gamma \rightarrow {\mathfrak g}_0$ denote the
  monodromy representation associated with the flat connection
  $\nabla$. Then:
  \begin{enumerate}\itemsep=0pt
  \item[$(2)$] 
  $\Gamma$ acts on ${\mathfrak g}_0$ by Lie algebra
    automorphisms. In particular, the monodromy representation
    integrates to a group homomorphism $\hat\mu \colon \Gamma
    \rightarrow \automorphism(G_0)$, where $G_0$ denotes the
    simply-connected Lie group with Lie algebra ${\mathfrak g}_0 $.
  \end{enumerate}

  Now let ${\mathfrak h}_0 \subset {\mathfrak g}_0={\mathfrak g}|_{m_0}$
  denote the kernel of the restriction of the anchor $\# \colon
  \mathfrak g \rightarrow TM$ to~$\mathfrak g|_{m_0}$, and let $H_0
  \subset G_0$ denote the connected subgroup with Lie algebra
  ${\mathfrak h}_0$. Then:
  \begin{enumerate}\itemsep=0pt
  \item[$(3)$] 
  The connection $\nabla $ is geometrically closed if
    and only if $H_0 \subset G_0$ is closed, in which case there
    exists a smooth action
    \[
    (\phi,x)\mapsto \phi_{G_0/H_0}(x) \colon \
    \Gamma \times {G_0/H_0} \rightarrow {G_0/H_0}
    \] of $\Gamma $ on
    ${G_0/H_0} $ by affine transformations, such that $\phi_{G_0/H_0}
    $ has twist $\hat \mu_\phi$.
  \item[$(4)$]
   There is an atlas of charts on $M$ with model ${G_0/H_0}$ such
    that each transition function is a restriction of some
    $\phi_{G_0/H_0}$, $\phi \in \Gamma$.
  \end{enumerate}
\end{PropositionBlaom}

\begin{RemarkBlaom}\label{Remark5.2}
  In the proposition statement (but not the subsequent proof)
  monodromy is to be understood in the usual sense of parallel
  translation (or analytic continuation) around closed paths; see the
  Remark~\ref{Remark2.6}.
\end{RemarkBlaom}

Before turning to the proof, let us explain how to recover genuine
local homogeneity from the conclusion of the proposition and hence
complete the proof of Theorem~\ref{TheoremA}. Indeed, let $G=\Gamma
\times_{\hat\mu} G_0$ denote the semidirect product, with
multiplication given by \begin{equation*}
  (\phi_1,g_1)(\phi_2,g_2):=\big( \phi_1
  \phi_2,\hat{\mu}_{\phi_2}^{-1}(g_1)g_2 \big).
\end{equation*}
Then, as $\phi_{G_0/H_0} $ is $G_0$-equivariant with twist $\hat
\mu_\phi$, the group $G$ acts on ${G_0/H_0}$ according to
$(\phi,g_0)\cdot x = \phi_{G_0/H_0}(g_0 \cdot x)$ ($x \in
G_0/H_0$). As this action is also transitive, we have a canonical
isomorphism ${G_0/H_0} \cong G/H$, where $H$ is the isotropy at $\mathrm{id} H_0
\in {G_0/H_0} $ of the action by $G$ just def\/ined. Viewing each
transition function as a local transformation of $G/H$, it becomes a~left-translation by some element of $\Gamma \cong \Gamma \times
\{\identity\}\subset G$.

\subsection{Reconstruction}\label{hutty}

To prove the proposition, pull ${\mathfrak g} $ back to a transitive
Lie algebroid $\tilde {\mathfrak g} $ over the universal cover~$\tilde
M$ of $M$, and $\nabla $ back to a Cartan connection $\tilde \nabla $
on $\tilde {\mathfrak g}$, as described in Section~\ref{holiday}. Then, by
Proposition~\ref{Proposition2.5}, the f\/inite-dimensional Lie algebra
${\mathfrak g}_0 \subset \Gamma(\tilde {\mathfrak g})$ of $\tilde
\nabla $-parallel sections acts on $\tilde M$ and this action is
transitive. That this Lie algebra may be identif\/ied with the f\/ibre
${\mathfrak g}|_{m_0}$, equipped with the bracket described in Proposition~\ref{Proposition5.1}(1),
follows from Remark~\ref{Remark2.6} and Proposition~\ref{Proposition2.6}(4)
(but we make no further use of this
interpretation).

Identify $\Gamma$ with the group of covering transformations and let
us understand the monodromy representation $\phi \mapsto \mu_\phi
\colon \Gamma \rightarrow \automorphism ({\mathfrak g}_0)$ in the
invariant sense def\/ined in Section~\ref{holiday}. That Proposition~\ref{Proposition5.1}(2) 
holds is
just Proposition~\ref{Proposition2.6}(2).
According to Proposition~\ref{Proposition2.6}(3),
each transformation $\phi \in \Gamma $ is a
${\mathfrak g}_0$-equivariant dif\/feomorphism of $M$ with twist
$\mu_\phi$.

It is clear from Proposition~\ref{Proposition2.6}(4)
that geometric
closure, in the sense of Section~\ref{geomclose}, is equivalent to the
condition in~Proposition~\ref{Proposition5.1}(3) 
above. If $\nabla $ is geometrically closed,
we can def\/ine the development $D \colon \tilde M \rightarrow
{G_0/H_0}$ associated with the action of ${\mathfrak g}_0$ on $\tilde
M$, as determined by the choice of some f\/ixed point $\tilde m_0 \in
\tilde M$ covering $m_0$. Def\/ining $\phi_{G_0/H_0} $ as in~\eqref{abox},
 one obtains an action as described in~Proposition~\ref{Proposition5.1}(3) 
above, on account of Proposition~\ref{Proposition4.7}(1) and~(2).  

The development $D \colon \tilde M \rightarrow {G_0/H_0} $ is a local
dif\/feomorphism (Proposition~\ref{Proposition4.3}) 
and, according to
Proposition~\ref{Proposition4.7}(3), 
transforms according to
\begin{gather}
  D(\phi(\tilde m))=\phi_{G_0/H_0} (D(\tilde m)),
  \qquad \phi \in \Gamma,\qquad \tilde m \in \tilde M.  \label{cvg}
\end{gather}

Now cover $M$ with open sets $U_i$, $i \in I$, small enough that each
set $U_i$ is evenly covered by some open set $\tilde U_i \subset
\tilde M$, and such that the restriction of the development $D
\colon \tilde U_i \rightarrow {G_0/H_0}$ is a~dif\/feomorphism onto its
image. Ref\/ining this covering if necessary, we may arrange that each
non-empty intersection $U_i \cap U_j $ is connected (see, e.g., \cite[Theorem~I.5.1]{Bott_Tu_82}). Let $s^i \colon U_i \rightarrow \tilde
U_i$ denote the inverse of the restriction $\tilde U_i \rightarrow
U_i$ of the covering $\tilde M \rightarrow M$. Def\/ine charts $\psi^i
\colon U_i \rightarrow {G_0/H_0}$ by $\psi^i=D \circ s^i$. Then, whenever
$U_i \cap U_j\ne \varnothing $, there exists a covering transformation
$\phi^{ij}\in \Gamma $ such that $s^j=\phi^{ij} \circ s^i$ on $U_i\cap
U_j$. (Here one uses the fact that two local continuous sections of a
covering map that have a common {\em connected} domain will agree on
the entire domain if they agree at one point.) Furthermore, for any $
m \in U_i\cap U_j$, we have, by~\eqref{cvg},
\begin{gather*}
  \psi^j(m)=D\big(s^j(m)\big)=D\big(\phi^{ij}(s^i(m))\big)
  =\phi^{ij}_{G_0/H_0}\big(D(s^i(m))\big)=\phi^{ij}_{G_0/H_0}\big(\psi^i(m)\big).
\end{gather*}
Whence the maps $\psi^i \colon U_i \rightarrow {G_0/H_0}$, $i \in I$,
constitute an atlas of charts meeting the requirement in Proposition~\ref{Proposition5.1}(4) 
above, and this concludes the proof of
Proposition~\ref{Proposition5.1}.

\section{Completeness}\label{completenessSection}
\subsection{Completeness in terms of the associated Lie algebra action}\label{recharacterize}

Let $\mathfrak g$ be a Lie algebroid on $M$ (not necessarily
transitive), let $\nabla $ be a f\/lat Cartan connection on~$\mathfrak
g$, and consider the associated action of the Lie algebra ${\mathfrak
  g}_0 $ on $\tilde M$ discussed in Section~\ref{holiday}.
\begin{PropositionBlaom}\label{Proposition6.1}
  A flat Cartan connection $\nabla$ is complete, in the sense of
  Section~{\rm \ref{coco}}, if and only if ${\mathfrak g}_0 $ acts on $\tilde M$
  by complete vector fields.
\end{PropositionBlaom}

\begin{proof}
  Adopting the notation of Section~\ref{holiday}, it is easy to see that
  $\nabla $ is complete if and only if the pullback connection $\tilde
  \nabla $ on $\tilde {\mathfrak g} $ is complete. But, according to
  Proposition~\ref{Proposition2.5}(4),
  $\tilde {\mathfrak g} $ is
  isomorphic to the action algebroid ${\mathfrak g}_0 \times \tilde
  M$, the connection $\tilde \nabla $ being represented by the
  canonical f\/lat connection on ${\mathfrak g}_0 \times \tilde M$. But
  geodesics of the canonical f\/lat connection on ${\mathfrak g}_0
  \times \tilde M$ are evidently those paths of the form
  $X_t=(\xi,\tilde m_t)$, where $\xi \in {\mathfrak g}_0 $ and $\tilde
  m_t$ is an integral curve of the corresponding vector f\/ield
  $\xi^\dagger$.
\end{proof}

Compactness of $M$ is insuf\/f\/icient to guarantee completeness:
\begin{CounterexampleBlaom}
  Let $M=S^1$ be the circle and let ${\mathfrak g}_0={\mathbb R}$ act
  on $\tilde M= {\mathbb R} $ according to
  $1^\dagger=e^{-\theta}\frac{\partial }{\partial \theta}$. Here
  $\theta$ denotes the standard coordinate function on ${\mathbb
    R}$. Evidently, $1^\dagger$ is not a~complete vector f\/ield.

  Def\/ine a representation $\mu \colon \Gamma \rightarrow
  \automorphism({\mathfrak g}_0)\cong {\mathbb R}\backslash \{0\}$ of
  the group of covering transformations $\Gamma\cong {\mathbb Z}$ on
  ${\mathfrak g}_0$ by $\mu_n = e^{2\pi n}$.  As a covering
  transformation, each $n \in \Gamma$ is the map $\theta \mapsto
  \theta + 2 \pi n$, which is a ${\mathfrak g}_0$-equivariant map with
  twist~$\mu_n$. By the discussion at the end of Section~\ref{holiday}, the
  quotient ${\mathfrak g} = (\tilde M \times {\mathfrak g}_0)/\Gamma$
  is a Lie algebroid over $M=S^1$ supporting a f\/lat Cartan connection,
  but the corresponding Lie algebra action of~${\mathfrak g}_0$ on~$\tilde M$ is, by construction, incomplete.
\end{CounterexampleBlaom}

\subsection{Suf\/f\/icient conditions for completeness}\label{suff}
Although compactness of $M$ is not suf\/f\/icient to guarantee
completeness, compactness {\em plus} simple-connectivity is obviously
suf\/f\/icient, for then $\tilde M$ is also compact (and all vector f\/ields
on compact manifolds are complete). More generally, we have:

 \begin{PropositionBlaom}\label{Proposition6.2}
   Let ${\mathfrak g} $ be a Lie algebroid over $M$ equipped with a
   flat Cartan connection $\nabla$, and let $\mu \colon \Gamma
   \rightarrow \automorphism({\mathfrak g}_0)$ denote the associated
   monodromy representation, as described in Section~{\rm \ref{holiday}}. Assume
   that one of the following conditions holds:
 \begin{enumerate}\itemsep=0pt
 \item[$(1)$] $M$ is compact and the image $\mu(\Gamma)\subset
   \automorphism({\mathfrak g}_0)$ of the monodromy representation has
   compact closure, denoted $\bar \Gamma \subset
   \automorphism({\mathfrak g}_0)$; or 
 \item[$(2)$] $M$ admits a ${\mathfrak g}$-invariant Riemannian metric
   $\sigma$ and this metric makes $M$ into a complete metric space
   $($automatically true if~$M$ is compact$)$.
 \end{enumerate}
 Then $\nabla $ is complete.
\end{PropositionBlaom}

Recall that because we assume $\nabla $ is Cartan, the
associated ${\mathfrak g}$-connection $\bar\nabla$ on $TM$ def\/ined in
Section~\ref{skoj} is in fact a $\mathfrak g$-representation. There is a
corresponding representation of ${\mathfrak g} $ on $S^2(T^*\!M)$, of
which $\sigma$ is a section. This explains what is meant by $\mathfrak
g$-invariance in~(2). Note also that~(1) is automatic in
the special case that $\automorphism({\mathfrak g}_0)$ is already a
compact group.

Before proving the preceding proposition, note that by Proposition~\ref{Proposition6.1} 
and Theorem~\ref{Theorem2.6}, it suf\/f\/ices to
establish the completeness of any action $\xi \mapsto \xi^\dagger$ of
some Lie algebra ${\mathfrak g}_0 $ on $\tilde M$, under the assumption
that each covering transformation $\phi \in \Gamma$ is ${\mathfrak
  g}_0$-equivariant with twist $\mu_\phi$. In this case the
${\mathfrak g}$-invariance of the metric $\sigma $ in hypothesis~(2) 
amounts to ${\mathfrak g}_0$-invariance of the lifted
metric on $\tilde M$.

\begin{proof}[Proof that (1) implies completeness]
  Fix some $\xi\in {\mathfrak g}_0$. We want to show that
  $\xi^\dagger$ is a~complete vector f\/ield on $\tilde M$. We will do
  this by f\/inding local lower bounds on the time of validity of its
  f\/low that are uniform with respect to the action of $\Gamma$; the
  compactness of~$M$ will then imply these bounds are globally
  uniform. To begin with, let us record that
  \begin{gather}
    \phi_*\xi^\dagger = (\mu_\phi \xi)^\dagger,\label{sore}
  \end{gather}
  for all $\xi \in {\mathfrak g}_0$ and $\phi \in \Gamma$; here
  $\phi_*$ denotes pushforward. This is just a restatement of the
  ${\mathfrak g}_0$-equivariance condition on the group of covering
  transformations.

  By the hypothesis (1), there exists a $\Gamma$-invariant
  inner product on ${\mathfrak g}_0$ (take an arbitrary inner product
  and average with respect to the Haar measure on the compact
  topological group $\bar \Gamma$). Equipping $M$ with an arbitrary
  Riemannian metric, we evidently obtain a $\Gamma$-invariant metric
  on~$\tilde M$.  At each $\tilde m\in \tilde M$ we now def\/ine
  \begin{gather*}
    c(\tilde m):=\sup_{\xi\neq 0}\frac{|\xi^\dagger(\tilde m)|}{|\xi |},
  \end{gather*}
  where the norms in the numerator and denominator are def\/ined with
  respect to the metric and inner product just introduced. It is clear
  $c(\tilde m)$ is a continuous with respect to~$\tilde m$.

With the help of~\eqref{sore}, $\Gamma$-invariance of the metric on
  $\tilde M$, as well as $\Gamma$-invariance of the inner product on
  ${\mathfrak g}_0$, one next shows that $c$ is $\Gamma$-invariant. We
  can therefore view it as a continuous function $c\colon M \rightarrow
  [0,\infty)$, and we have the $\Gamma$-uniform estimate
\begin{gather}
  |\xi^\dagger(\tilde m)|\le c(\pi(\tilde m))|\xi|,
  \qquad\xi\in {\mathfrak g}_0,\qquad \tilde m\in \tilde M. \label{esti}
\end{gather}
Here $\pi\colon \tilde M \rightarrow M$ denotes the covering
projection.

We will now use the following lemma, whose proof is an elementary
exercise in Riemannian geometry.
  \begin{LemmaBlaom}\label{Lemma6.4}
    Let $m$ be an arbitrary point of a Riemann manifold $ M$, and let
    $r>0$ be small enough that the open geodesic ball~$B_{3r}(m)$ of
    radius $3r$ about $m$ is well-defined\footnote{Here we mean a
      geodesic=normal ball in the sense of, e.g., p.~70 of~\cite{doCarmo_92}.}. Then for any vector field $V$ on $M$,
    integral curves of $V$ beginning in $B_{r}(m)$ are well-defined
    for all times $t\in[-T,T]$, where $T:=r/\|V\|_{m,2r}$. Here,
    $\|V\|_{m,2r}$ denotes the supremum of $|V(n)|$ over all $n\in
    B_{2r}(m)$.
  \end{LemmaBlaom}

  So, let $m\in M$ be arbitrary. Then there is $r=r(m)>0$ such that
  the geodesic ball $B_{3r}(m)\subset M$ is well-def\/ined, has compact
  closure, and is evenly covered by a disjoint union of geodesic balls
  $B_{3r}(\tilde m) \subset \tilde M$, one for each $\tilde m$ lying
  over $m$ in the covering. For any such $\tilde m$ we have, in the
  notation of the lemma,
  \begin{gather*}
    \|\xi^\dagger\|_{\tilde m,2r}\le C(m)|\xi|,\qquad\text{where}\quad C(m):=\sup\{c(m)\suchthat m\in B_{2r}(m)\}.
  \end{gather*}
  Here we have applied~\eqref{esti} and have $C(m)<\infty$ because $c$
  is continuous and $B_{2r}(m)$ has compact closure. Applying the
  lemma to the manifold $\tilde M$, we conclude that for any initial
  condition covering a point in $B_r(m)\subset M$, the corresponding
  integral curve of $\xi^\dagger$ is def\/ined for all times $t\in
  [-T(m),T(m)]$, where $T(m):=r(m)/(C(m)|\xi|)>0$. Since $M$ is
  compact, there is a f\/inite covering $\{B_{r(m_j)}(m_j)\}_{j=1}^N$ of
  $M$, for some $m_1,\ldots, m_N\in M$. So an integral curve of
  $\xi^\dagger$ with {\em any} initial condition is def\/ined for all
  times $t\in[-T_\mathrm{min},T_\mathrm{min}]$, where
  $T_\mathrm{min}:=\min_{j=1}^N T(m_j)>0$. Since this lower bound on
  times is independent of the initial condition, the curve is def\/ined
  for all time.
\end{proof}

\begin{proof}[Proof that (2) implies completeness]
  In this case every vector f\/ield $\xi^\dagger$ is a Killing f\/ield on
the universal cover $\tilde M$, to which the metric $\sigma $
  lifts. In particular, every integral curve $t \mapsto \tilde m(t)
  \colon (a,b) \rightarrow \tilde M $ of $\xi^\dagger$ has constant
  speed. Since $M$ is a complete metric space, so is~$\tilde M$.

  Suppose that $t_n \rightarrow b<\infty$, for some sequence $t_1,
  t_2, \ldots \in (a,b)$. Because $m(t)$ has constant speed, $m(t_1),
  m(t_2), \ldots$ is a Cauchy sequence, which must therefore converge
  to some $m_b\in M$. Any integral curve of $\xi^\dagger$ beginning at
  $m_b$ may, by uniqueness of integral curves, be regarded as an
  extension of $m(t)$, i.e., $m(t)$ extends to some time interval
  $(a,b_+)$ with $b_+>b$. A similar argument applies to $a$. Whence
  $\xi^\dagger$ is complete by the usual proof by contradiction.
\end{proof}

\subsection{The proof of Theorem~\ref{TheoremB}}\label{completetheorem}

With Proposition~\ref{Proposition6.1} in hand, the necessity of
completeness in Theorem~\ref{TheoremB} is not dif\/f\/icult to prove, and we now turn
to the proof of suf\/f\/iciency. To that end, we begin with a variant of
the result sought in which the model space $G_0/H_0$ is the quotient
of a simply-connected Lie group by a connected subgroup:
\begin{PropositionBlaom}\label{Proposition6.5}
  Let $\mathfrak g$ be a transitive Lie algebroid, over a smooth
  connected manifold~$M$, and~$\nabla $ a~flat, geometrically closed,
  complete, Cartan connection on~$\mathfrak g$. Then:
  \begin{enumerate}\itemsep=0pt
  \item[$(1)$]
  The subgroup $H_0 \subset G_0$ is closed, when $G_0$ and $H_0$ have the meanings given in Proposition~{\rm \ref{Proposition5.1}}, and the universal cover of~$M$ is~$G_0/H_0$.

  \item[$(2)$] 
  The group $\Gamma \cong \pi_1(M)$ of covering
    transformations is a subgroup of $\affine({G_0/H_0})$, and each
    element $\phi \in \Gamma $ has twist $\hat \mu_\phi$. Here $\phi
    \mapsto \hat \mu_\phi \colon \Gamma \rightarrow
    \automorphism(G_0)$ is a group homomorphism determined by the
    monodromy of $\nabla $ and defined in Proposition~{\rm \ref{Proposition5.1}}.
  \end{enumerate}
\end{PropositionBlaom}

To recover Theorem~\ref{TheoremB} from this variant, i.e., to obtain $M \cong
(G/H)/\Gamma $, with $\Gamma \subset G$ acting by left translations,
one def\/ines $G$ to be the semidirect product $\Gamma \times_{\hat \mu}
G_0 $ as described already after the statement of Proposition~\ref{Proposition5.1}.

\begin{proof}
  Given a pair $({\mathfrak g}, \nabla)$ satisfying the hypotheses of
  the proposition, consider the correspon\-ding action of the Lie
  algebra ${\mathfrak g}_0$ on the universal cover $\tilde M$ detailed
  in Section~\ref{holiday}. In particular, $G_0$ is then the
  simply-connected Lie group integrating ${\mathfrak g}_0$. According
  to Proposition~\ref{Proposition6.1} and Palais' integrability
  theorem~\cite{Palais_57a}, the action of ${\mathfrak g}_0$ can be
  integrated to an action of the Lie group $G_0$. Therefore, in the
  def\/inition of the development $D \colon \tilde M \rightarrow
  G_0/H_0$ of the ${\mathfrak g}_0$-action~-- see Section~\ref{devmap}, but
  read $\tilde M$ in place of~$M$~-- we may take ${\mathcal G}$ to be
  the action groupoid $G_0 \times \tilde M$, and the groupoid morphism
  $\Omega \colon {\mathcal G} \rightarrow G_0 $ is just the projection
  $\Omega \colon G_0 \times \tilde M \rightarrow G_0$.
  The subgroup $H_0 \subset G_0$ is closed because it is the isotropy subgroup of the $G_0$-action at~$\tilde m_0$.
   Moreover,
  $P=G_0 \times \{\tilde m_0\}$ and the development $D \colon \tilde M
  \rightarrow G_0/H_0$ is seen to be a dif\/feomorphism. This
  establishes~(1).

  In the subsequent identif\/ication $\tilde M \cong {G_0/H_0}$, each
  covering transformation $\phi \in \Gamma$ and corresponding map
  $\phi_{G_0/H_0} $ appearing in Proposition~\ref{Proposition5.1} is
  identif\/ied. The former therefore have the properties claimed in~(2).
\end{proof}

\section{Illustrations}\label{illustration}
Following are simple illustrations of the theory expounded in the
present article.

\subsection{Local Lie groups}\label{when}
As in \cite{Blaom_05}, we call pair of linear connections $\nabla $,
$\bar\nabla $ on a smooth manifold~$M$ (more precisely, on its tangent
bundle $TM$) a {\it dual pair} if
$\bar\nabla_XY-\nabla_YX=[X,Y]$. Notice that a connection $\nabla $ is
a~Cartan connection precisely when its dual $\bar\nabla $ is f\/lat. If
both~$\nabla $ and~$\bar\nabla $ are f\/lat, it follows from Theorem~\ref{TheoremA}
that~$M$ is a locally homogeneous manifold modelled on some Lie group~$G$. For the purposes of this section, we accordingly def\/ine a {\it
  local Lie group} to be a smooth connected manifold~$M$, equipped with
a dual pair of simultaneously f\/lat connections~$\nabla $,~$\bar\nabla
$.

A bone f\/ide Lie group $G_0$ is a local Lie group in this sense: the
canonical f\/lat connec\-tions~$\nabla $,~$\bar\nabla $~-- corresponding to right and left trivialization, respectively~--
furnish the required dual pair. In this case, both connections $\nabla
$ and $\bar\nabla $ are complete (the geodesics being the right and
left cosets of one-parameter subgroups of $G$) and have trivial
monodromy. Conversely, we have:

\begin{PropositionBlaom}[cf.~\S~3.8.3 of~\cite{Sharpe_97}]\label{Proposition7.1}
  Let $(M, \nabla, \bar\nabla)$ be a local Lie group. Then:
  \begin{enumerate}\itemsep=0pt
  \item[$(1)$]
  The torsion, $\torsion\bar\nabla=-\torsion \nabla$ of
    $\bar\nabla $, is a $\bar\nabla$-parallel tensor whose restriction
    to any tangent space ${\mathfrak g}_0 = T_{m_0}M$, $m_0 \in M$, is
    a Lie bracket.
  \end{enumerate}

  If, additionally, $\nabla $ is complete, and $G_0$ denotes the
  simply-connected Lie group integrating~${\mathfrak g}_0$, then:
\begin{enumerate}\itemsep=0pt
  \item[$(2)$]
  The universal cover of $M$ is $G_0$ and the group of covering
    transformations $\Gamma \cong \pi_1(M)$ acts by affine
    transformations.
  \end{enumerate}

  If $\nabla $ is both complete and has trivial monodromy, then:
\begin{enumerate}\itemsep=0pt
  \item[$(3)$] 
  There exists an embedding of $\Gamma $ into $G_0$ such that
    $M \cong G_0/\Gamma$.
  \end{enumerate}
\end{PropositionBlaom}
\begin{proof}
  The claim (1) is just a special case of Proposition~\ref{Proposition5.1}(1).
  Conclusion~(2) is a special case
  of Proposition~\ref{Proposition6.5}. 
  Under the additional
  hypothesis that $\nabla $ has trivial monodromy,
 Proposition~\ref{Proposition6.5}(2) 
 implies that each element of $\Gamma
  \subset \affine(G_0)$ has trivial twist. Each such element is
  therefore a left translation, allowing us to identify $\Gamma $ with
  a subgroup of $G_0$, as claimed in~(3).
\end{proof}

Of course, in the case that the embedding $\Gamma \subset G_0$ in~(3) 
is normal, $M$ is a global Lie group.

The problem of globalizing a local Lie group structure is also
considered in \cite{Abadoglu_Ortacgil_10} (under the additional
assumption that the connection $\bar\nabla $ comes from a {\it global}
parallelism on $M$). In particular, it is shown that a necessary
condition for globalizability is the vanishing of the class $[w] \in
H^1(M)$, where $w$ is the closed one-form def\/ined by
$w(U)=\trace(\torsion \bar\nabla (U,\,\cdot \,))$, which is closed on
account of Bianchi's second identity and the hypothesis $\curvature
\nabla =0$.

\subsection{Riemannian manifolds with constant scalar curvature}
Let $M$ be a smooth connected manifold. Associated with any Riemannian
metric $\sigma$ on $M$ is a natural Lie subalgebroid ${\mathfrak g}
\subset J^1 ({TM})$ which, according to \cite{Blaom_12}, supports a
canonical Cartan connection $\nabla $. The curvature of $\nabla $
vanishes on any simply-connected region $U$ precisely when the vector
space of Killing f\/ields on $U$ is of maximal possible dimension, and
the Lie algebra of all such Killing f\/ields can be concretely
described. We now give a global analogue of this result.

By def\/inition, ${\mathfrak g} $ is the subbundle of all one-jets
$J^1_mV$ of vector f\/ields $V$ on $M$ such that $\sigma $ has vanishing
Lie derivative at the point $m \in M$. To describe the Cartan
connection $\nabla $ explicitly, let $\levi$ denote the Levi-Cevita
connection, which we may regard as a splitting of the canonical exact
sequence,
\begin{gather*} 0\rightarrow T^*\!M\otimes TM\rightarrow
  J^1 (TM)\rightarrow TM\rightarrow 0.
\end{gather*}
There is a corresponding exact sequence,
\begin{gather*}
  0\longrightarrow  {\mathfrak h} \longrightarrow {\mathfrak g} \longrightarrow  TM \longrightarrow 0,
\end{gather*}
where ${\mathfrak h}\subset T^*\!M\otimes TM$ denotes the
$\mathfrak{o}(n)$-bundle of all $\sigma$-skew-symmetric tangent space
endomorphisms, which is $\levi$-invariant. Because $\sigma $ is
$\levi$-invariant, the splitting subbundle of $J^1(TM)$ lies inside
$\mathfrak g$ and we obtain an identif\/ication ${\mathfrak g}\cong
TM\oplus{\mathfrak h}$. Under this identif\/ication, the Cartan
connection $\nabla $ is given by
\begin{gather*}
  \nabla_U(V\oplus\phi)=\big(\levi_UV+\phi (U)\big)\oplus
  \big(\levi_U\phi+\curvature\levi  (U, V)\big).
\end{gather*}
With the help of the Bianchi identities for linear connections, one
computes
\begin{gather*}
        \curvature\nabla (U_1, U_2)(V\oplus\phi)=0\oplus
        \big({-}\big(\levi_V\curvature\levi+\phi\cdot\curvature\levi\big)(U_1,U_2)\big),
\end{gather*}
implying implying that $\nabla$ is f\/lat if and only if $\curvature\levi$
is both $\levi$-invariant and ${\mathfrak h}$-invariant.  Now
${\mathfrak h} $-invariance implies, by purely algebraic arguments,
that $\curvature\levi$ has only a scalar component:
\begin{gather}
        \curvature\levi (V_1, V_2)=s\big( \sigma (V_2)\otimes V_1-\sigma
        (V_1)\otimes V_2 \big),\label{pop}
\end{gather}
for some function $ s\in{\mathbb R}$, the {\it scalar curvature};
$\levi$-invariance implies $s$ is a constant.

With a view to applying Proposition~\ref{Proposition6.5} 
in the f\/lat
case, we f\/irst def\/ine ${\mathfrak g}_0$ and ${\mathfrak h}_0 \subset
{\mathfrak g}_0 $ as in Proposition~\ref{Proposition5.1}, and, to obtain
the bracket on ${\mathfrak g}_0$, compute the torsion, $\torsion
\bar\nabla $, of the associated ${\mathfrak g}$-connection
$\bar\nabla$ on ${\mathfrak g} $ (which is actually a ${\mathfrak
  g}$-representation). To this end, one requires a formula for the Lie
bracket on ${\mathfrak g} \subset J^1(TM)$. Indeed, under the identif\/ication
${\mathfrak g} \cong {TM} \oplus {\mathfrak h} $, we have
\begin{gather}
  [V_1\oplus\phi_1,V_2\oplus\phi_2]_{\mathfrak g}=
  [V_1, V_2]_{{TM}}\oplus \big( [\phi_1,\phi_2]_{\mathfrak
    h}+\levi_{V_1}\phi_2-\levi_{V_2}\phi_1 +
  \curvature\levi(V_1,V_2) \big)\label{piu},
\end{gather} where
$[\,\cdot\,,\,\cdot \,]_{\mathfrak h}$ is the bracket on the
$\mathfrak{o}(n)$-bundle ${\mathfrak h} $. This formula is an instance
of the general formula \cite[Proposition~6.2(4)]{Blaom_12}. Supposing
that $\nabla $ is f\/lat, and hence~\eqref{pop} holds, one obtains with
the help of~\eqref{piu},
\begin{gather*}
  \torsion\bar{\nabla} (V_1\oplus\phi_1,
  V_2\oplus\phi_2)=
  ( \phi_1(V_2)-\phi_2(V_1) ) \oplus
  \big( [\phi_1,\phi_2]_{\mathfrak h}-s(\sigma (V_2)\otimes
  V_1+\sigma (V_1)\otimes V_2) \big).
\end{gather*}

Restricting $\torsion \bar\nabla $ to obtain the Lie bracket on
${\mathfrak g}_0={\mathfrak g}|_{m_0}$, we see from the formula that
${\mathfrak g}_0$ is the Lie algebra of $\mathfrak{o}(n)\times
{\mathbb R}^n $ (semidirect product), $\mathfrak{o}(n+1)$, or
$\mathfrak{o}(n,1)$, according to whether $s=0$, $s>0$, or $s<0$,
respectively. In each case ${\mathfrak h}_0 \subset {\mathfrak g}_0$
is isomorphic to $\mathfrak{o}(n)$.

By def\/inition, the metric $\sigma $ is ${\mathfrak g}$-invariant, so
that completeness of $M$, as a metric space, implies completeness of
the Cartan connection $\nabla $, by Proposition~\ref{Proposition6.2}. In that
case, Proposition~\ref{Proposition6.5} 
establishes the following:

\begin{PropositionBlaom}[cf.~Proposition~8.4.3 of~\cite{doCarmo_92}]\label{Proposition7.2}
  If a connected Riemannian manifold $M$ is a~complete metric space
  and has constant scalar curvature $s$. Then:
  \begin{enumerate}\itemsep=0pt
    \item[$(1)$] The universal cover of $M$ is:
      \begin{itemize}\itemsep=0pt
      \item $G_0/H_0\cong {\mathbb R}^n$, in the case $s=0$, where $G_0$
        is the simply-connected Lie group with Lie algebra ${\mathbb
          R}^n \times \mathfrak{o}(n)$ $($semidirect product$)$;
      \item $G_0/H_0 \cong S^n$, in the case $s>0$, where $G_0$ is the simply-connected
        Lie group with Lie algebra $\mathfrak{o}(n+1)$;
      \item $G_0/H_0\cong {\mathbb H}^n$, in case $s<0$, where $G_0$
        is the simply-connected Lie group with Lie algebra
        $\mathfrak{o}(n,1)$.
      \end{itemize}
      In every case $H_0 \subset G_0$ is the connected subgroup with
      Lie algebra $\mathfrak{o}(n)$, realised as a~subalgebra of
      ${\mathfrak g}_0$ in the usual way.
    \item[$(2)$] The group of covering transformations $\Gamma \cong
      \pi_1(M)$ acts on $G_0/H_0$ by affine transformations.
    \end{enumerate}
\end{PropositionBlaom}

\vspace{-2.5mm}

\subsection*{Acknowledgements} The author acknowledges many helpful
discussions with Erc\"ument Orta\c{c}gil.
\vspace{-2mm}

\pdfbookmark[1]{References}{ref}
\LastPageEnding


\begin{thebibliography}{99}
\footnotesize\itemsep=-1pt

\bibitem{Abadoglu_Ortacgil_10}
Abado{\u{g}}lu E., Orta{\c{c}}gil E., Intrinsic characteristic classes of a
  local {L}ie group, \href{http://dx.doi.org/10.4171/PM/1873}{\textit{Port. Math.}} \textbf{67} (2010), 453--483,
  \href{http://arxiv.org/abs/0912.0855}{arXiv:0912.0855}.

\bibitem{Alekseevsky_Michor_95b}
Alekseevsky D.V., Michor P.W., Dif\/ferential geometry of {$\mathfrak
  g$}-manifolds, \href{http://dx.doi.org/10.1016/0926-2245(95)00023-2}{\textit{Differential Geom. Appl.}} \textbf{5} (1995), 371--403,
  \href{http://arxiv.org/abs/math.DG/9309214}{math.DG/9309214}.

\bibitem{Blaom_05}
Blaom A.D., Geometric structures as deformed inf\/initesimal symmetries,
  \href{http://dx.doi.org/10.1090/S0002-9947-06-04057-8}{\textit{Trans. Amer. Math. Soc.}} \textbf{358} (2006), 3651--3671,
  \href{http://arxiv.org/abs/math.DG/0404313}{math.DG/0404313}.

\bibitem{Blaom_12}
Blaom A.D., Lie algebroids and {C}artan's method of equivalence, \href{http://dx.doi.org/10.1090/S0002-9947-2012-05441-9}{\textit{Trans.
  Amer. Math. Soc.}} \textbf{364} (2012), 3071--3135, \href{http://arxiv.org/abs/math.DG/0509071}{math.DG/0509071}.

\bibitem{Blaom_C}
Blaom A.D., A Lie theory of pseudogroups, {i}n preparation.

\bibitem{Bott_Tu_82}
Bott R., Tu L.W., Dif\/ferential forms in algebraic topology, \textit{Graduate
  Texts in Mathematics}, Vol.~82, Springer-Verlag, New York, 1982.

\bibitem{Crainic_Fernandes_11}
Crainic M., Fernandes R.L., Lectures on integrability of {L}ie brackets, in
  Lectures on {P}oisson Geometry, \textit{Geom. Topol. Monogr.}, Vol.~17, Geom.
  Topol. Publ., Coventry, 2011, 1--107, \href{http://arxiv.org/abs/math.DG/0611259}{math.DG/0611259}.

\bibitem{Dazord_97}
Dazord P., Groupo\"\i de d'holonomie et g\'eom\'etrie globale, \href{http://dx.doi.org/10.1016/S0764-4442(97)80107-3}{\textit{C.~R.
  Acad. Sci. Paris S\'er.~I Math.}} \textbf{324} (1997), 77--80.

\bibitem{doCarmo_92}
do~Carmo M.P., Riemannian geometry, \textit{Mathematics: Theory {\rm \&} Applications},
  Birkh\"auser Boston Inc., Boston, MA, 1992.

\bibitem{Goldman_10}
Goldman W.M., Locally homogeneous geometric manifolds, in Proceedings of the
  {I}nternational {C}ongress of {M}athematicians, {V}ol.~{II}, Hindustan Book
  Agency, New Delhi, 2010, 717--744.

\bibitem{Gunning_67}
Gunning R.C., Special coordinate coverings of {R}iemann surfaces, \href{http://dx.doi.org/10.1007/BF01362287}{\textit{Math.
  Ann.}} \textbf{170} (1967), 67--86.

\bibitem{Kamber_Michor_04}
Kamber F.W., Michor P.W., Completing {L}ie algebra actions to {L}ie group
  actions, \href{http://dx.doi.org/10.1090/S1079-6762-04-00124-6}{\textit{Electron. Res. Announc. Amer. Math. Soc.}} \textbf{10}
  (2004), 1--10, \href{http://arxiv.org/abs/math.DG/0310308}{math.DG/0310308}.

\bibitem{Kowalski_97}
Kowalski O., On strictly locally homogeneous {R}iemannian manifolds,
  \href{http://dx.doi.org/10.1016/S0926-2245(96)00043-5}{\textit{Differential Geom. Appl.}} \textbf{7} (1997), 131--137.

\bibitem{Mackenzie_05}
Mackenzie K.C.H., General theory of {L}ie groupoids and {L}ie algebroids,
  \textit{London Mathematical Society Lecture Note Series}, Vol.~213, Cambridge
  University Press, Cambridge, 2005.

\bibitem{Palais_57a}
Palais R.S., A global formulation of the {L}ie theory of transformation groups,
  \textit{Mem. Amer. Math. Soc.} \textbf{22} (1957), iii+123~pages.

\bibitem{Sharpe_97}
Sharpe R.W., Dif\/ferential geometry. Cartan's generalization of Klein's Erlangen
  program, \textit{Graduate Texts in Mathematics}, Vol.~166, Springer-Verlag,
  New York, 1997.

\bibitem{Thurston_97}
Thurston W.P., Three-dimensional geometry and topology, {V}ol.~1,
  \textit{Princeton Mathematical Series}, Vol.~35, Princeton University Press,
  Princeton, NJ, 1997.

\end{thebibliography}
\end{document}